\documentclass{article}

\usepackage{graphicx}
\usepackage[fleqn]{amsmath}
\usepackage{amsthm}
\usepackage{amssymb}
\usepackage{amsbsy}
\usepackage{amsfonts}
\usepackage{amstext}
\usepackage{amscd}
\usepackage[dvips]{epsfig}
\allowdisplaybreaks

\ifx\undefined\bysame
\newcommand{\bysame}{\leavevmode\hbox to3em{\hrulefill}\,}
\fi

\newcommand{\al}{\alpha}

\newcommand{\bp}[2]{\begin{picture}(#1)#2\end{picture}}
\newcommand{\C}{{\Bbb C}}
\newcommand{\cA}{{\cal A}}
\newcommand{\cB}{{\cal A}(\ast)}

\newcommand{\cP}{{\cal P}}
\newcommand{\cQ}{{\check Q}}
\newcommand{\cZ}{{\check Z}}

\newcommand{\Dg}{\Delta_{\fg^*}}
\newcommand{\Dh}{\Delta_{\fh^*}}
\newcommand{\e}{\varepsilon}

\newcommand{\ff}{{\bold f}}
\newcommand{\fg}{{\frak g}}
\newcommand{\fh}{{\frak h}}
\newcommand{\fns}{\footnotesize}
\newcommand{\fsl}[1]{{\frak s}{\frak l}_{#1}}
\newcommand{\gC}{{{\frak g}_{{}_\C}}}
\newcommand{\hC}{{{\frak h}_{{}_\C}}}

\newcommand{\hLMO}{\hat Z^{\mbox{\scriptsize L{\hspace{-1pt}}M{\hspace{-1pt}}O}}}

\newcommand{\hWg}{{\hat W}_\fg}
\newcommand{\JdA}[2]{\begin{picture}(30,10)\put(0,-1){\pc{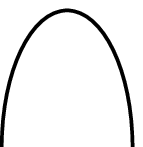}{0.5}}
  \put(3,-10){\scriptsize $#1$}\put(21,-10){\scriptsize $#2$}\end{picture}}
\newcommand{\JdU}[2]{\begin{picture}(32,10)\put(0,-4){\pc{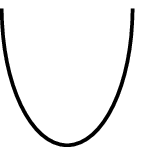}{0.5}}
  \put(1.5,13.5){\scriptsize $#1$}\put(20,13.5){\scriptsize $#2$}\end{picture}}

\newcommand{\la}{\lambda}
\newcommand{\La}{\Lambda}
\newcommand{\lra}{\longrightarrow}
\newcommand{\ns}{\normalsize}

\newcommand{\ot}{\otimes}
\newcommand{\ptl}[1]{\partial_{#1}}
\newcommand{\R}{{\Bbb R}}
\newcommand{\Ups}{\Upsilon}

\newcommand{\vspc}[1]{\begin{picture}(0,#1)\end{picture}}

\newcommand{\Wg}{W_\fg}
\newcommand{\Z}{{\Bbb Z}}

\newcommand{\pc}[2]{\mbox{$\begin{array}{c}
   \includegraphics[scale=#2]{figure/#1.eps}
   \end{array}$}}
\newcommand{\pcp}[4]{\begin{picture}(#1)\put(#2){$\pc{#3}{#4}$}\end{picture}}

\theoremstyle{plain}
   \newtheorem{thm}{Theorem}[section]
   \newtheorem{lem}[thm]{Lemma}
   \newtheorem{prop}[thm]{Proposition}
   
\theoremstyle{definition}

   \newtheorem{rem}[thm]{Remark}

   \newcommand{\sn}{\operatorname{sign}}
   \newcommand{\vol}{\operatorname{Vol}}
   \newcommand{\ad}{\operatorname{ad}}
  \newcommand{\cM}{{\mathcal C}}
  \newcommand{\cI}{{\mathcal I}}
  \newcommand{\cF}{{\mathcal F}}
  \newcommand{\cH}{{\mathcal H}}
  \newcommand{\Gr}{\mathrm{Gr}}
  \newcommand{\od}{{\mathcal D}}
  \newcommand{\Tr}{{\mathrm{Tr}}}
  \newcommand{\qd}{\mbox{{\normalsize q}-}\!{\cal D}}
  \newcommand{\BQ}{\mathbb{Q}}
  \newcommand{\BR}{\mathbb{R}}
  \newcommand{\eL}{\mathcal{E}}
  \newcommand{\teL}{\tilde {\mathcal{E}}}
  \newcommand{\cO}{{\mathcal O}}
  \newcommand{\psdraw}[2]
         {\begin{array}{c} \hspace{-1.3mm}
         \raisebox{-4pt}{\psfig{figure=#1.eps,width=#2}}
         \hspace{-1.9mm}\end{array}}

\begin{document}

\vspace*{2pc}

{\bf\Large \begin{center}
The perturbative invariants of rational homology 3-spheres \\
can be recovered from the LMO invariant
\end{center}}

\vskip 1.5pc

\centerline{\Large Takahito Kuriya, \ Thang T. Q. Le, \ Tomotada Ohtsuki}

\large

\vskip 2pc

\begin{abstract}
We show that
the perturbative $\fg$ invariant of rational homology 3-spheres
can be recovered from the LMO invariant for any simple Lie algebra $\fg$,
{\it i.e.},
the LMO invariant is universal among the perturbative invariants.
This universality was conjectured in \cite{LMO}.
Since the perturbative invariants dominate the quantum invariants
of integral homology 3-spheres
\cite{Habiro_q,Habiro_WRT,HabiroLe},
this implies that
the LMO invariant dominates the quantum invariants
of integral homology 3-spheres.
\end{abstract}

\vskip 1pc

\section{Introduction}

\begin{figure}[htpb]
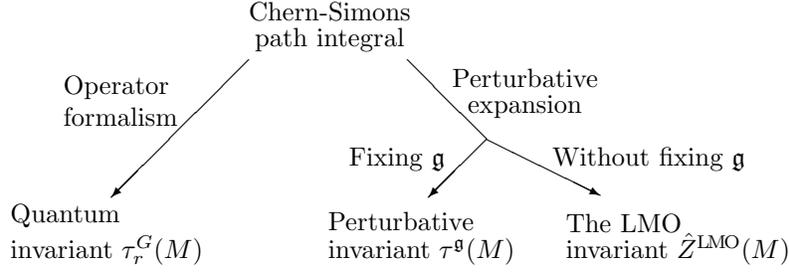

$$
\bp{300,100}{
\put(100,80){\shortstack[c]{\ns Chern-Simons \\ \ns path integral}}
\put(100,75){\vector(-1,-1){52}}
\put(30,50){\shortstack[l]{\ns Operator \\ \ns formalism}}
\put(10,0){\shortstack[l]{\ns Quantum \\ \ns invariant $\tau^G_r(M)$}}
\put(190,45){\line(-1,1){30}}
\put(177,55){\shortstack[c]{\ns Perturbative \\ \ns expansion}}
\put(190,45){\vector(-1,-1){22}}
\put(138,35){\ns Fixing $\fg$}
\put(215,35){\ns Without fixing $\fg$}
\put(130,0){\shortstack[l]{\ns Perturbative \\ \ns invariant $\tau^\fg(M)$}}
\put(190,45){\vector(2,-1){43}}
\put(220,0){\shortstack[l]
   {\ns The LMO \\ \ns invariant \bp{10,0}{\put(0,0){\ns $\hLMO(M)$}}}}}
$$
\caption{\label{fig.pb} Physical background}
\end{figure}

\begin{figure}[htpb]
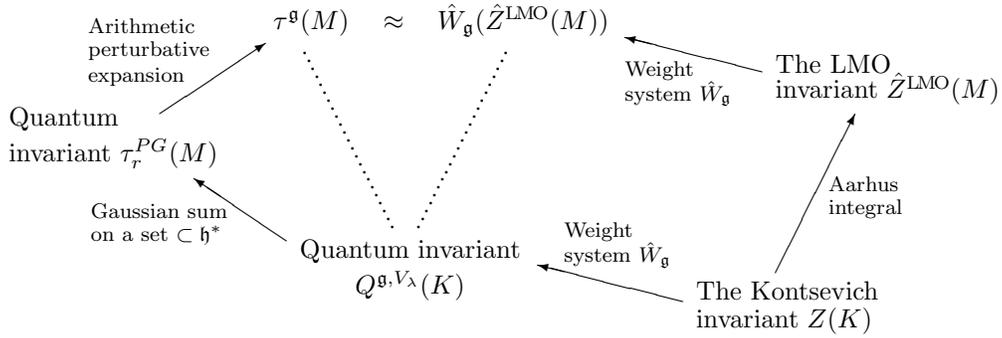

$$
\bp{320,120}{
\put(80,100){\ns $\tau^\fg(M) \quad\approx\quad \hWg(\hLMO(M))$}
\put(37,68){\vector(3,2){40}}
\put(10,80){\shortstack[l]
   {\fns Arithmetic \\ \fns perturbative \\ \fns expansion}}
\put(-20,50){\shortstack[l]{\ns Quantum \\ \ns invariant $\tau^{P G}_r(M)$}}
\put(265,84){\vector(-4,1){52}}
\put(213,73){\shortstack[l]
   {\fns Weight \\ \fns system \bp{10,0}{\put(0,0){\fns $\hWg$}}}}
\put(270,74){\shortstack[l]
   {\ns The LMO \\ \ns invariant \bp{10,0}{\put(0,0){\ns $\hLMO(M)$}}}}
\put(85,20){\vector(-3,2){35}}
\put(11,20){\shortstack[l]{\fns Gaussian sum \\ \fns on a set $\subset \fh^*$}}
\multiput(125,25)(-1.5,3){23}{\circle*{1}}
\multiput(135,25)(1.5,3){23}{\circle*{1}}
\put(90,0){\shortstack[c]{\ns Quantum invariant \\ $Q^{\fg,V_\la}(K)$}}
\put(235,-3){\vector(-4,1){55}}
\put(190,12){\shortstack[l]
   {\fns Weight \\ \fns system \bp{10,0}{\put(0,0){\fns $\hWg$}}}}
\put(270,8){\vector(1,2){30}}
\put(290,30){\shortstack[l]{\fns Aarhus \\ \fns integral}}
\put(240,-13){\shortstack[l]{\ns The Kontsevich \\ \ns invariant $Z(K)$}}}
$$
\caption{\label{fig.mc} Mathematical construction}
\end{figure}

In the late 1980s,
Witten \cite{Witten} proposed
topological invariants of a closed 3-manifold $M$
for a simple compact Lie group $G$,
which is formally presented by
a path integral whose Lagrangian is the Chern-Simons functional
of $G$ connections on $M$.
There are two approaches to obtain mathematically rigorous information
from a path integral:
the operator formalism and the perturbative expansion.
Motivated by the operator formalism of the Chern-Simons path integral,
Reshetikhin and Turaev \cite{RT}
gave the first rigorous mathematical construction of quantum invariants
of 3-manifolds,
and, after that, rigorous constructions of quantum invariants
of 3-manifolds were obtained by various approaches.
When $M$ is obtained from $S^3$ by surgery along a framed knot $K$,
the quantum $G$ invariant $\tau^G_r(M)$ of $M$
is defined to be a linear sum of
the quantum $(\fg,V_\la)$ invariant $Q^{\fg,V_\la}(K)$ of $K$
at an $r$th root of unity,
where $\fg$ is the Lie algebra of $G$,
and $V_\la$ denotes the irreducible representation of $\fg$
whose highest weight is $\la$.
On the other hand,
the perturbative expansion of the Chern-Simons path integral suggests that
we can obtain the perturbative $\fg$ invariant
(a power series) when we fix $\fg$,
and obtain the LMO invariant (an infinite linear sum of trivalent graphs)
when we make the perturbative expansion without fixing $\fg$.
As a mathematical construction,
we can define the perturbative $\fg$ invariant $\tau^\fg(M)$
of a rational homology 3-sphere $M$
by arithmetic perturbative expansion of $\tau_r^{PG}(M)$
as $r \to \infty$
\cite{Ohtsuki_pi,Rozansky,Le_isp},
where $PG$ denotes the quotient of $G$ by its center.
Further, we can present the LMO invariant $\hLMO(M)$ \cite{LMO}
of a rational homology 3-sphere $M$
by the Aarhus integral \cite{Aint3}.
It was conjectured \cite{LMO} that
the perturbative $\fg$ invariant can be recovered from the LMO invariant
by the weight system $\hWg$
for any simple Lie algebra $\fg$.
In the $\fsl2$ case, this has been shown in \cite{Ohtsuki_rec}.
See Figures \ref{fig.pb} and \ref{fig.mc},
for these invariants and relations among them.

The aim of this paper is to show the following theorem.

\begin{thm}[{\rm see \cite{Aint1,Kuriya}}\footnote{
It was announced in \cite{Aint1} that
the perturbative $\fg$ invariant can be recovered from the LMO invariant.
However, their proof is not published yet.
The first author \cite{Kuriya} showed a proof,
but his proof is partially incomplete.
The aim of this paper is to show a complete proof of the theorem.
}]
\label{thm.rec}
Let $\fg$ be any simple Lie algebra.
Then, for any rational homology 3-sphere $M$,
$$
\hWg \big( \hLMO (M) \big)
\ = \
| H_1(M;\Z) |^{({\rm dim} \, \fg \, - \, {\rm rank} \, \fg)/2} \, \tau^\fg(M) ,
$$
where $| H_1(M;\Z) |$ denotes the cardinality of
the first homology group $H_1(M;\Z)$ of $M$.
\end{thm}

\noindent
We give two proofs of the theorem:
a geometric proof (Sections \ref{sec.geom_pf_knot} and \ref{sec.link_fti})
and an algebraic proof (Sections \ref{sec.alg_pf_knot} and \ref{sec.link_alg}).
The theorem implies that
the LMO invariant dominates the perturbative invariants.
Further, since the perturbative invariants dominate the quantum Witten-Reshetikhin-Turaev invariants
of integral homology 3-spheres
\cite{Habiro_q,Habiro_WRT,HabiroLe},
it follows from the theorem that
the LMO invariant dominates the quantum invariants
of integral homology 3-spheres.\footnote{
For rational homology 3-spheres,
it is known \cite{BBL} that
the quantum WRT invariant $\tau_r^{SO(3)}(M)$, at roots of unity of order co-prime to the order of the first homology group,
can be obtained from the perturbative invariant $\tau^{\fsl2}(M)$.
Hence, the LMO invariant $\hLMO(M)$ dominates
$\tau_r^{SO(3)}(M)$ for those roots of unity.}

\begin{figure}[htpb]
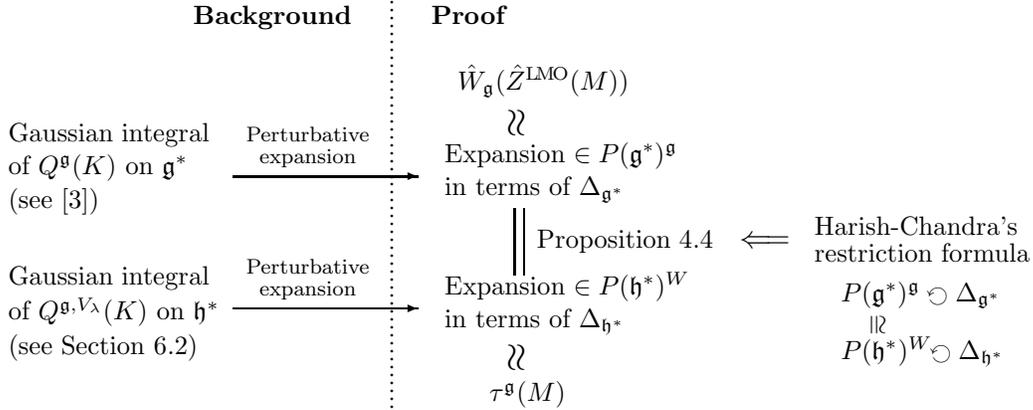

$$
\bp{360,140}{
\put(155,110){\ns $\hWg(\hLMO(M))$}
\put(174,92){\Large \rotatebox{90}{$\approx$}}
\put(167,-8){\ns $\tau^\fg(M)$}
\put(174,3){\Large \rotatebox{90}{$\approx$}}
\multiput(130,-10)(0,3){52}{\circle*{1}}
\put(55,135){\ns\bf Background}
\put(145,135){\ns\bf Proof}
\put(-15,65){\shortstack[l]{\ns Gaussian integral \\
  \ns of $Q^{\fg}(K)$ on $\fg^*$ \\ \ns (see \cite{Aint1})}}
\put(-15,10){\shortstack[l]{\ns Gaussian integral \\
  \ns of $Q^{\fg,V_\la}(K)$ on $\fh^*$ \\ \ns (see Section \ref{sec.p_pi})}}
\put(70,77){\vector(1,0){70}}
\put(70,27){\vector(1,0){70}}
\put(75,82){\shortstack[c]{\fns Perturbative \\ \fns expansion}}
\put(75,32){\shortstack[c]{\fns Perturbative \\ \fns expansion}}
\put(150,70){\shortstack[l]
   {\ns Expansion $\in P(\fg^*)^\fg$ \\ \ns in terms of $\Dg$}}
\put(150,20){\shortstack[l]
   {\ns Expansion $\in P(\fh^*)^W$ \\ \ns in terms of $\Dh$}}
\put(176.5,40){\line(0,1){25}}
\put(180,40){\line(0,1){25}}
\put(185,50){\ns Proposition \ref{prop.DhD} \ \ \large $\Longleftarrow$}
\put(290,45){\shortstack[l]
   {\ns Harish-Chandra's \\ \ns restriction formula}}
\put(300,30){\ns $P(\fg^*)^\fg \
\raisebox{-2pt}[0pt][0pt]{\rotatebox{90}{$\circlearrowleft$}} \ \Dg$}
\put(310,25){\ns \rotatebox{-90}{$\cong$}}
\put(300,8){\ns $P(\fh^*)^W
\raisebox{-2pt}[0pt][0pt]{\rotatebox{90}{$\circlearrowleft$}} \ \Dh$}}
$$
\caption{\label{fig.s_ap}
Sketch of the algebraic proof of Theorem \ref{thm.rec},
when $M$ is obtained from $S^3$ by surgery along a framed knot $K$}
\end{figure}

Let us  explain a sketch of the  proof when $M$ is obtained by surgery on a knot.
The LMO invariant $\hLMO (M)$ can be presented by
the Aarhus integral \cite{Aint3}.
It is shown from this presentation that
the image $\hWg \big( \hLMO (M) \big)$ can be presented
by an integral of Gauss type over the dual $\fg^*$, or alternatively by an expansion
given in terms of the Laplacian $\Dg$ of $\fg^*$.
On the other hand, as we explain in Section \ref{sec.p_pi}, the perturbative invariant
$\tau^\fg(M)$ is presented by a Gaussian integral over $\fh^*$,
where $\fh$ is a Cartan subalgebra of $\fg$, or alternatively by  an expansion
given in terms of the Laplacian $\Dh$ of $\fh^*$.
We then show that $\hWg \big( \hLMO (M) \big)= \tau^\fg(M)$ by establishing a result relating integrals  over $\fg^*$
and integrals over $\fh^*$, similar to the well known Weyl reduction integration formula. Alternatively, we  show  $\hWg \big( \hLMO (M) \big)= \tau^\fg(M)$ by
using Harish-Chandra restriction theorem that relates the Laplacian $\Dg$ on $\fg^*$ to the Laplacian $\Dh$ on $\fh^*$.
For a sketch of the algebraic proof, see also Figure \ref{fig.s_ap}.

In case when $M$ is obtained by surgery on a link we also present two proofs. The first one is more algebraic. We reduce the theorem to the case of surgery on knots
by using the fact that the operators involved are invariant under the action of $\fg$. The other proof has  quite a different flavor. We show that two multiplicative finite type invariants of rational
homology spheres are the same if they agree on the set of rational homology spheres obtained by surgery on knots (for finer results see Theorem \ref{022}). This result is also interesting by itself.
The theorem then follows, since both $\hWg \big( \hLMO (M) \big)$ and $\tau^\fg(M)$, up to any degree, are finite type. This part relates the paper to the origin of the theory: The discovery of the perturbative invariant of homology 3-spheres for $SO(3)$ case \cite{Ohtsuki_pi} leads the third author to define finite type invariants of 3-manifolds.

The paper is organized, as follows.
In Section \ref{sec.def_p},
we review definitions of terminologies,
and show some properties of Jacobi diagrams.
In Section \ref{proof2},
we present the proof of the main theorem, based on results proved later.
We consider the knot case in Section \ref{sec.knot} and the link case in
 Section \ref{sec.link}.
In Section \ref{sec.ppi}, we discuss how the perturbative invariant can be obtained as an asymptotic expansion of the Witten-Reshetikhin-Turaev invariant,
and give a proof that our formula of the perturbative invariant is coincident with that given in \cite{Le_isp}. We also show that finite parts of the perturbative invariant $\tau^\fg$
are of finite type.

\vskip 1pc

The third author would like to thank Susumu Ariki
for pointing out Harish-Chandra's restriction formula
when he tried to prove Proposition \ref{prop.DhD} in the $\fsl3$ case.
The authors would like to thank
Dror Bar-Natan, Kazuo Habiro, Andrew Kricker,
Lev Rozansky, Toshie Takata and Dylan Thurston
for valuable comments and suggestions
on early versions of the paper.

\section{Preliminaries}
\label{sec.def_p}

In this section, we recall basic facts about Lie algebras
in Section \ref{0001}  and
theory of the Kontsevich invariant in Section \ref{sec.dp}.
We show some facts about Laplacian operators in Euclidean spaces
in Section \ref{Laplace},
and present the LMO invariant in Section \ref{sec.pr_LMO}.

\subsection{Lie algebra}
\label{0001}
In this paper, $G$ is a compact connected simple Lie group,
$\fg$ its Lie algebra and  $\fh$ a fixed Cartan  subalgebra of $\fg$.
Up to scalar multiplication,
there is a unique Ad-invariant inner product on $\fg$.
The complexification $\gC$ of $\fg$ can be presented as
$\gC= \fg + \sqrt{-1} \, \fg$.
Then $\hC= \fh + \sqrt{-1} \, \fh$ is a Cartan subalgebra of $\gC$.

There is a root system $\Phi_\C \subset (\hC)^*$
of the pair $(\gC,\hC)$.
It is known that $\Phi_\C$, as well as the weights of $\fg$-modules,
are purely imaginary, {\it i.e.},
$\Phi_\C \subset ( \sqrt{-1} \, \fh)^*\equiv \sqrt{-1} \, \fh^*$.
Following the common convention in Lie algebra theory
(see {\it e.g.} \cite{Hall}),
we call $\beta \in \fh^*$ a {\em real root}
(resp. a {\em real weight of a $\fg$-module})
if $\sqrt{-1} \, \beta$ is a root (resp. a weight of the $\fg$-module).
We normalize the invariant inner product
so that the square length of every short root is 2.
We denote by $W$ the Weyl group,
$\Phi_+$  the set of positive real roots of $\fg$, and
$\rho$ the half-sum of positive real roots. Let
$\phi_+$ be the number of positive roots of $\fg$. One has
$\phi _+ = ({\rm dim}\, \fg - {\rm dim}\, \fh)/2
= ({\rm dim}\, \fg - {\rm rank}\, \fg)/2$. We denote by
$V_\la$ the irreducible representation of $\fg$
whose highest weight is $\sqrt{-1} \, \la$.

Let $S(\fg)$ and  $U(\fg)$ be respectively the symmetric tensor algebra
and the universal enveloping algebra of $\fg$.
One can naturally identify $S(\fg)$ with
$P(\fg^*)$, the algebra of polynomial functions on $\fg^*$. Throughout the paper, $\hbar$ is a formal parameter, and $q=e^\hbar \in \R[[\hbar]]$.
One considers  $S(\fg)[[\hbar]]$ as a ring of functions on $\fg^*$ with values in $\R[[\hbar]]$.

The following $W$-skew-invariant functions $\od$ is important to us:
$$
\od(\lambda) \, := \, \prod_{\al \in \Phi_+} \frac{(\la,\al)}{(\rho,\al)} \, .
$$
When $\lambda-\rho$ is a dominant real weight,
$\od(\lambda)$ is
the dimension  of $V_{\la-\rho}$.

We identify $\fh^*$ with a subspace of $\fg^*$ using the
invariant inner product. For a function $g$ on $\fg^*$, its restriction to $\fh^*$
will be denoted by $\cP(g)$.

A source of function on $\fh^*$ is given by the enveloping algebra $U(\fg)$. For $g\in U(g)$ we define a polynomial function, also denoted by $g$, on $\fh^*$ as follows.
Suppose $\lambda-\rho$ is a dominant real weight. One can take the trace $\Tr_{V_{\lambda-\rho}}(g)$ of the action of $g$ in the $\fg$-module $V_{\lambda-\rho}$.
It is known that there is a unique polynomial function, denoted by also by $g$,   on $\fh^*$ such that $g(\lambda)=\Tr_{V_{\lambda-\rho}}(g)$.

There is a vector space isomorphism $\Upsilon_\fg: S(\fg) \to U(\fg)$, known as the Duflo-Kirillov map (see \cite{BGV,BGRT_w,BLT}). We can extend $\Upsilon_\fg$ multi-linearly to  a vector space isomorphism
$\Upsilon_\fg: S(\fg)^{\otimes \ell } \to U(\fg)^{\otimes \ell }$.
When restricted to the $\fg$-invariant parts, $\Upsilon_g:
S(\fg)^\fg \to U(\fg)^\fg$ is an algebra isomorphism. Note that $U(\fg)^\fg$ is the center of the algebra $U(\fg)$.

\subsection{Laplacian and Gaussian integral on a Euclidean space}
\label{Laplace}

 Let $V$ be a Euclidean space. In our applications we will always have $V=\fg$ or $V=\fh$ with the Euclidean structure coming from the
invariant inner product.
As usual, one identifies the symmetric algebra $S(V)$ with the polynomial function algebra $P(V^*)$.
 The Laplacian $\Delta_V$, associated with the Euclidean structure of $V$,  acts on $S(V)= P(V^*)$  and is defined by
$$
\Delta_V \, = \, \sum_{i} \ptl{x_i}^2,
$$
where $x_i$'s are coordinate functions with respect to
an orthonormal basis of $V$. It is known that for $x,y \in V$,  $\frac{1}{2}\Delta_V(xy)= (x,y)$, the inner product of $x$ and $y$.

Let $\hbar$ be a formal parameter. For a non-zero real number $f$ let us consider the following operator $\eL^{(f)}_V: S(V)= P(V^*) \to \R[1/\hbar]$ expressed through an exponent of the Laplacian and the evaluation at $0$:

$$ \eL^{(f)}_V(g) = \exp\left(-\frac{\Delta}{2f\hbar}\right)\, (g) \big|_{x=0} \, \in \R[1/\hbar].$$
Because $\Delta_V$ is a second order differential operator,
it is easy to see that if $g$ is a homogeneous polynomial, then
\begin{equation} \hspace*{8pc}
\eL^{(f)}_V(g) \ = \ \begin{cases}   0 &\text{ if $\deg(g)$ is odd,} \\
\frac{\text{scalar}}{\hbar ^{\deg(g)/2}} &\text{ if $\deg(g)$ is even.}
\end{cases}
\label{n100}
\end{equation}

Adjoining $\hbar$, we get an extension $\eL^{(f)}_V: S(V)(\!(\hbar)\!)= P(V^*)(\!(\hbar)\!) \to \R(\!(\hbar)\!)$ as follows. If $g = \sum_{n=-\infty }^\infty g_n \, \hbar^n$ with $g_n \in S(V)$, then
$$ \eL^{(f)}_V(g) = \sum_n  \eL^{(f)}_V (g_n) \hbar^n   \in \BR(\!(\hbar)\!).$$

There is a generalization to the multi-variable case.
Suppose $\ff := (f_1,\dots,f_\ell)$ is an $\ell$-tuple of non-zero real numbers and $g_1\otimes \dots \otimes g_\ell\in S(V)^{\otimes \ell}$, then we define

$$ \eL^{(\ff)}_V(g_1 \otimes \dots \otimes g_\ell) = \prod_{j=1}^\ell \eL^{(f_j)}_V (g_j).$$

Formally we can put $\eL^{(\ff)}_V = \bigotimes_j \eL^{(f_j)}_V$. Again there is an obvious extension $\eL^{(\ff)}_V: S(V)^{\otimes \ell}(\!(\hbar)\!) \to \R(\!(\hbar)\!)$.

\subsection{Jacobi diagrams, weight systems, and the Kontsevich invariant}
\label{sec.dp}

In this section,
we review  Jacobi diagrams, weight systems, and the Kontsevich invariant of framed string links.
For details see {\it e.g.} \cite{Ohtsuki_book}.

A {\it uni-trivalent graph}
is a graph every vertex of which is either univalent
or trivalent. A uni-trivalent graph is {\it vertex-oriented} if at
 each trivalent
 vertex a cyclic order of edges is fixed. For a 1-manifold $Y$,
a {\it Jacobi diagram} on $Y$ is
the manifold $Y$ together with
a vertex oriented uni-trivalent graph such that
univalent vertices of the graph are distinct points on $Y$.
In figures we draw $Y$ by thick lines and the uni-trivalent graphs
by thin lines, in such a way that
each trivalent vertex is vertex-oriented in the counterclockwise order.
We define the {\it degree} of a Jacobi diagram
to be half the number of univalent and trivalent vertices
of the uni-trivalent graph of the Jacobi diagram.
We denote by $\cA(Y)$
the quotient vector space
spanned by Jacobi diagrams on $Y$
subject to the following relations,
called {\it the AS, IHX, and STU relations} respectively,
$$
\pc{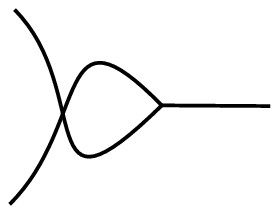}{0.4} = - \pc{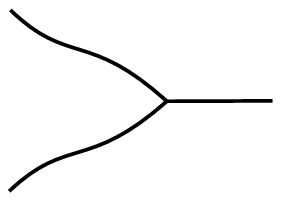}{0.4} , \quad
\pc{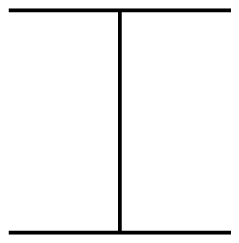}{0.4} = \pc{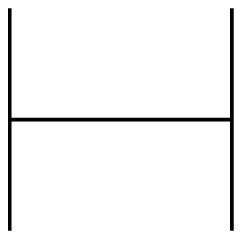}{0.4} - \pc{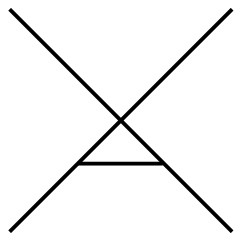}{0.4} , \quad
\pc{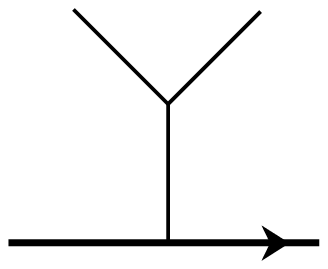}{0.4} = \pc{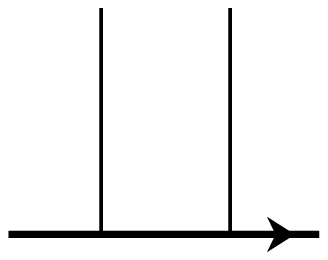}{0.4} - \pc{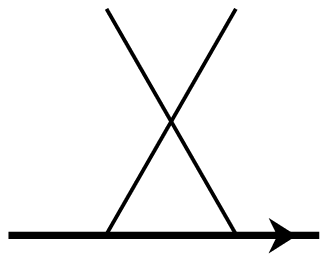}{0.4} .
$$
For $S= \{ x_1, \cdots, x_\ell \}$,
a {\it Jacobi diagram} on $S$ is a vertex-oriented uni-trivalent graph
whose univalent vertices are labeled by elements of $S$.
We denote by $\cA(\ast_S)$
the quotient vector space spanned by Jacobi diagrams on $S$
subject to the AS and IHX relations.
In particular, when $S$ consists of a single element,
we denote $\cA(\ast_S)$ by $\cA(\ast)$.
\ $\cA(\emptyset)$ and $\cA(\ast_S)$ form algebras
with respect to the disjoint union of Jacobi diagrams, and
$\cA( \sqcup^\ell \! \downarrow )$ forms an algebra
with respect to the vertical composition of copies of
$\sqcup^\ell \! \downarrow$.

We briefly review weight systems;
for details, see \cite{BarNatan,Ohtsuki_book}.
We define the {\it weight system} $\Wg(D)$ of a Jacobi diagram $D$
by ``substituting'' $\fg$ into $D$, {\it i.e.},
putting $D$ in a plane, $\Wg(D)$ is defined to be the composition
of intertwiners, each of which is given at each local part of $D$
as follows.
$$
\pcp{50,10}{0,-6}{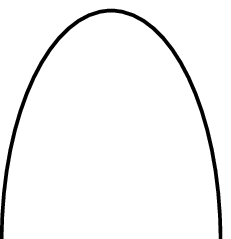}{0.5}
\begin{CD} \R \\ @AA{B}A \\ \fg \ot \fg \end{CD} \qquad\qquad
\pcp{50,10}{0,1}{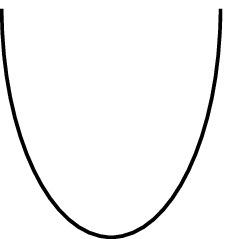}{0.5}
\begin{CD} \fg \ot \fg \\ @AAA \\ \R \end{CD} \qquad\qquad
\pcp{50,10}{0,-3}{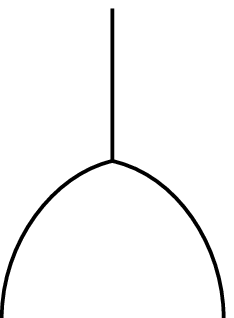}{0.5}
\begin{CD} \fg \\ @AA{[\,\cdot\,,\,\cdot\,]}A \\ \fg \ot \fg \end{CD}
$$
Here, the first map is the invariant form of $\fg$, and
the second map is the map taking $1$ to $\sum_i X_i \ot X_i$,
where $\{ X_i \}_{i \in I}$ is
an orthonormal basis of $\fg$ with respect to the invariant form,
and the third map is the Lie bracket of $\fg$.
For $D_1 \in \cA(\ast)$ and $D_2 \in \cA(\downarrow)$,
we have the following intertwiners
as the compositions of the above maps,
and we can define
$\Wg(D_1) \in S(\fg)$ and $\Wg(D_2) \in U(\fg)$
as the images of $1$ by these maps.
$$
\bp{60,10}{\put(0,0){\pc{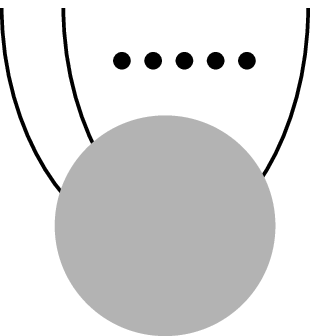}{0.5}}\put(22,-6){$D_1$}}
\begin{CD} S(\fg) \\ @AAA \\ \R \end{CD} \qquad\qquad
\bp{70,10}{\put(0,2){\pc{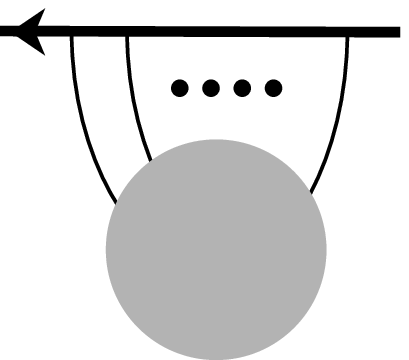}{0.5}}\put(29,-6){$D_2$}}
\begin{CD} U(\fg) \\ @AAA \\ \R \end{CD}
$$
In a similar way, we can also define
$\Wg : \cA(\sqcup^\ell \! \downarrow ) \to \big( U(\fg)^{\ot \ell} \big)^\fg$
and
$\Wg : \cA(\ast_S) \to \big( S(\fg)^{\ot \ell} \big)^\fg$;
they are algebra homomorphisms. Note that there is a standard degree on the polynomial algebra $S(\fg)^{\ot \ell}$ which carries over to $U(\fg)^{\ot \ell}$ by the Poincare-Birkhoff-Witt isomorphism.
If $D$ is a diagram with $k$ univalent vertices,
then $\Wg(D)$ has degree $\le k$.
The weight system $W_{\gC}$ is defined in the same way.
Since $W_{\gC} = \Wg$ by definition, we denote $W_{\gC}$ by $\Wg$.
Further, we define $\hWg$ by
$\hWg(D) = \Wg(D) \, \hbar^d$
for a Jacobi diagram $D$ of degree $d$.

There is a formal Duflo-Kirillov algebra isomorphism $ \Upsilon: \cA(\ast) \to \cA(\downarrow)$ (see \cite{BGRT_w,BLT}). The  obvious multi-linearly extension $\Upsilon: \cA(\{x_1,\dots,x_\ell\}) \to \cA(\sqcup^\ell\downarrow) $ is not an algebra isomorphism, but a vector space isomorphism.
The following diagram is
commutative \cite[Theorem 3]{BGRT_w}.

\begin{equation} \hspace*{7pc}
\label{eq.cd}
\begin{CD}
\cB @>{\hWg}>> S(\fg)^\fg[[\hbar]]
\begin{picture}(0,0) \put(13,5){\line(1,0){20}}\put(13,2){\line(1,0){20}}
\end{picture}
@.  P(\fg^*)^\fg[[\hbar]]
\\
@V{\Ups}V{\cong}V  @V{\Ups_\fg}V{\cong}V  @V{\cong}V{\cP}V
\\
\cA(\downarrow) @>>{\hWg}> U(\fg)^\fg[[\hbar]] @>{\cong}>{\psi_\fg}>
P(\fh^*)^W[[\hbar]]
\end{CD}
\end{equation}
Here, $P(\fh^*)^W$ denotes
the algebra of $W$-invariant polynomial functions on $\fh^*$.
$\Upsilon$ denotes the Duflo-Kirillov isomorphism.
$\psi_\fg$ denotes the Harish-Chandra isomorphism;
for $\la \in \fh^*$,
$\psi_\fg(z)(\la)$ is defined to be the scalar by which
$z \in U(\fg)^\fg$ acts on the irreducible representation of $\fg$
whose highest weight is $\la-\rho$. In other words, $\psi_\fg(z) (\lambda) = z(\lambda)/\od(\lambda)$.
$\cP$ is the restriction map from $\fg^*$ to $\fh^*$.

A {\it string link} is an embedding $\varphi$ of
$\ell$ copies of the unit interval,
$[0,1] \times \{ 1 \}$, $\cdots$, $[0,1] \times \{ \ell \}$,
into $[0,1] \times \C$,
so that $\varphi \big( (\e,j) \big) = (\e,j)$
for all $\e \in \{ 0,1 \}$ and $1 \le j \le \ell$.
We obtain a link from a string link
by closing each component of $\sqcup^\ell \! \downarrow$.
A (string) link is called {\it algebraically split}
if the linking number of each pair of components is $0$.

The {\it Kontsevich invariant} $Z(T)$ \cite{Kontsevich,LM}
of an $\ell$-component framed string link $T$
is defined to be in $\cA(\sqcup^\ell \! \downarrow)$;
for its construction, see, {\it e.g.}, \cite{LM,Ohtsuki_book}.
Let $\nu = Z(U)$, the Kontsevich invariant of the unknot $U$ with framing 0; the exact value of $\nu$ is calculated in \cite{BLT}.
Using the Poincare-Birkhoff-Witt isomorphism $\cA(S^1) \cong \cA(\downarrow)$
(see \cite{BarNatan}),
we will consider $\nu$ as an element in $\cA(\downarrow)$.

Let $\Delta^{(\ell)}: \cA(\downarrow) \to \cA(\sqcup^\ell \! \downarrow)$
be the cabling operation which replaces an arrow by $\ell$ parallel copies
(see {\it e.g.} \cite[Section 1]{LM}).
The modification $\cZ(T)$ of $Z(T)$ used in the definition of
the LMO invariant is
$$
\cZ(T) := \nu^{\otimes \ell} \, \left ( \Delta^\ell (\nu) \right) \, Z(T).
$$
Applying $\Ups^{-1}$ followed by the weight map, we define the following element:

\begin{equation} \hspace*{8pc}
\cQ^\fg (T)= \hWg\Big(\Ups^{-1} \big(\cZ(T)\big) \Big) \in  \big( S(\fg)^{\ot \ell} \big)^\fg [[\hbar]].
\label{n2000}
\end{equation}

\subsection{Presentations of the LMO invariant}
\label{sec.pr_LMO}

In this section,
we recall and modify a formula of the LMO invariant \cite{LMO} of a rational homology 3-sphere $M$
using the Aarhus integral \cite{Aint3}  for the case when $M$ is obtained by surgery along an algebraically split link.

Suppose $T$ is an algebraically split $\ell$-component string link with 0 framing on each component, and $L$ is its closure. Suppose the components of $T$ are ordered. Let $\ff=(f_1,\dots,f_\ell)$ be an $\ell$-tuple of
$\ell$ non-zero integers, and $M$ be the rational homology 3-sphere obtained by surgery on $L$ with framing $f_1,\dots, f_\ell$.

Let $\theta \in \cA(\emptyset)$ be the following Jacobi diagram
\begin{equation} \hspace*{12pc}
\label{eq.theta}
\theta \ = \, \pc{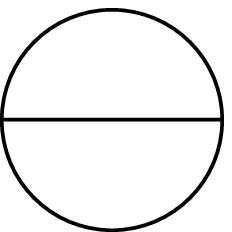}{0.5} \in \cA(\emptyset) .
\end{equation}
Define
\begin{equation} \cI(T,\ff):= \exp \Big( - \frac{\sum_j f_j}{48} \theta \Big)
\Bigg\langle \prod_j \exp \Big( -\frac1{2 f_j}
\JdU{\ptl{x_j}}{\ptl{x_j}} \Big), \,
\big( \Ups^{-1} \big) \cZ(T) \Bigg\rangle
\in \cA(\emptyset) .
\label{defn_I}
\end{equation}
Here, for a Jacobi diagram $D_1$ whose univalent vertices are
labeled by $\ptl{x_j}$'s
and a Jacobi diagram $D_2$ whose univalent vertices are
labeled by $x_j$'s,
we define the bracket by
$$
\langle D_1, D_2 \rangle
\, = \,
\Big(
\begin{array}{l}
\mbox{the sum of all ways of gluing the $\ptl{x_j}$-labeled univalent vertices}
\\
\mbox{of $D_1$ to the $x_j$-labeled univalent vertices of $D_2$ for each $j$}
\end{array} \Big) \in \cA(\emptyset) ,
$$
if the number of $\ptl{x_j}$-labeled univalent vertices of $D_1$
are equal to the number of $x_j$-labeled univalent vertices of $D_2$
for each $j$, and put $\langle D_1, D_2 \rangle = 0$ otherwise.
In particular, when $T= \downarrow$ is the trivial string link, one has
$$
\cI(\downarrow,\pm 1) \, = \,
\exp \Big( \mp \frac{1}{48} \theta \Big)
\Bigg\langle \exp \Big( \mp \frac1{2} \JdU{\ptl{x}}{\ptl{x}} \Big), \,
\Ups^{-1}(\nu^2) \Bigg\rangle \ \in \cA(\emptyset) ,
$$
Then, the LMO invariant of $M$ is presented by\footnote{
The bracket of this presentation is called the Aarhus ``integral'',
since its corresponding Lie algebra version is actually an integral
on $(\fg^*)^{\oplus \ell}$ \cite{Aint1}.}
\begin{equation} \hspace*{8pc}
\label{eq.p_LMO}
\hLMO(M) \, = \, \frac{ \cI(T,\ff)}{  \prod_{j=1}^\ell \cI(\downarrow,\sn(f_j))}
\in \cA(\emptyset) .
\end{equation}

We remark that
the presentation (\ref{eq.p_LMO}) is obtained from \cite[Theorem 6]{BaLa},
noting that (with notations from \cite{BaLa})
\begin{align*}
& {\rm \AA}_0(L) \, = \,
\int \Big( \prod_j \Ups_{x_j}^{-1} \Big)\big( \cZ(L) \big) \, d X ,
\\
& \Big( \prod_j \Ups_{x_j}^{-1} \Big)\big( \cZ(L) \big)
\, = \,
\Big( \prod_j \Ups_{x_j}^{-1} \Big)\big( \cZ(T) \big)
\exp \Big( - \frac{\sum_j f_j}{48} \theta \Big)
\prod_j \exp \Big( \frac{f_j}{2} \JdA{x_j}{x_j} \Big),
\end{align*}
which are obtained from
Lemma 3.8 and Corollaries 3.11 and 3.12 of \cite{BaLa}.

\section{Proof of the main theorem} In this section we show the proof of the main theorem in Section \ref{sec.proof}
based on results proved in later sections.
\label{proof2}

\subsection{Comparing the LMO invariant and the perturbative invariant} We again assume $M,L,T,\ff$ the same as in Section \ref{sec.pr_LMO}.
Recall that in \eqref{n2000} we defined $\cQ^\fg(T) \in (S(\fg)^{ \otimes \ell})^\fg[[\hbar]]$.

\begin{prop} \label{n2005}
Assume the above notations.
\\
\noindent
{\rm a)} The LMO invariant of $M$, after applied by the weight map, has the following presentation
\begin{equation} \hspace*{8pc}
\hWg(\hLMO(M)) = \frac{I_1(T,\ff)}{ {  \prod_{j=1}^\ell I_1(\downarrow,\sn(f_j))}}\, ,
\label{n2001}
\end{equation}
where
$$I_1(T,\ff) =  \left( \prod_{j=1}^\ell  q^{-f_j |\rho|^2/2} \right) \,
\eL^{(\ff)}_\fg \left ( \cQ^\fg(T) \right).$$
\\
\noindent
{\rm b)} The perturbative invariant has the following presentation
\begin{equation} \hspace*{10pc}
\tau^\fg(M) = \frac{I_2(T,\ff)}{ {  \prod_{j=1}^\ell I_2(\downarrow,\sn(f_j))}}\, ,
\label{n20011}
\end{equation}
where
$$
I_2(T,\ff) = \left( \prod_{j=1}^\ell  q^{-f_j |\rho|^2/2} \right) \,
\eL^{(\ff)}_\fh \left( \od^{\otimes \ell}\, \Upsilon_\fg ( \cQ^\fg(T)) \right).
$$

\end{prop}
\begin{proof}
Apply the algebra map $\hWg$ to \eqref{eq.p_LMO},
$$
\hWg(\hLMO(M))=
\frac{ \hWg (\cI(T,\ff)) }
{ \prod_{j=1}^\ell  \hWg \big ( \cI(\downarrow,\sn(f_j))\big) } \, .
$$
Using Lemmas \ref{lem.W_U}, \ref{Wg_theta} and the definition of $\cI(T,\ff)$ in \eqref{defn_I} we get
$$
\hWg (\cI(T,\ff)) = I_1(T,\ff),
$$
which proves part (a) of the proposition. Part (b) will be proved
in Section \ref{sec.p_Le}.
\end{proof}

To prove the main theorem one needs to understand the relation between $\eL^{(\ff)}_\fg$ and $\eL^{(\ff)}_\fh$. We will prove the following proposition
in Sections \ref{sec.knot} and \ref{sec.link}.
\begin{prop} There is a non-zero constant $c_\fg$ such that for  $g\in (S(\fg)^{\otimes \ell})^\fg[[\hbar]]$ and any $\ell$-tuple $\ff=(f_1,\dots,f_\ell)$ of non-zero integers one has
\begin{equation} \hspace*{6pc}
\eL^{(\ff)}_\fg \left(  g \right) = \left ( \prod_{j=1}^\ell    (-2f_j \hbar)^{\phi_+}\, c_\fg\right) \, \eL^{(\ff)}_\fh \left(   \od^{\otimes \ell}\, \Upsilon_\fg ( g) \right).
\label{n2008}
\end{equation}
\label{n2003}
\end{prop}

\subsection{Proof of Main Theorem} \label{sec.proof}
Now we can prove Theorem \ref{thm.rec}. First we assume that $M$ can be obtained by surgery along an algebraically split link $L$. We assume $T,\ff$
as in Section \ref{sec.pr_LMO}.
 One has
\begin{align} I_1(T,\ff) & =  \left( \prod_{j=1}^\ell  q^{-f_j |\rho|^2/2} \right) \, \eL^{(\ff)}_\fg \left ( \cQ^\fg(T) \right) \notag\\
&= \left( \prod_{j=1}^\ell  q^{-f_j |\rho|^2/2} \right) \, \left ( \prod_{j=1}^\ell    (-2f_j \hbar)^{\phi_+}\, c_\fg\right) \, \eL^{(\ff)}_\fh \left(   \od^{\otimes \ell}\, \Upsilon_\fg \left ( \cQ^\fg(T) \right) \right) \notag\\
&= \left ( \prod_{j=1}^\ell    (-2f_j \hbar)^{\phi_+}\, c_\fg\right) \, I_2(T,\ff), \label{n2007}
\end{align}
where the second equality follows from Proposition \ref{n2003} since $\cQ^\fg(T) \in (S(\fg)^{\otimes \ell})^\fg[[\hbar]]$.
In particular, applying \eqref{n2007} for $(T,\ff)=(\downarrow,\sn(f_j))$, then taking the product when $j$ runs from 1 to $\ell$, one has
\begin{equation} \hspace*{3pc}
\prod_{j=1}^\ell I_1(\downarrow,\sn(f_j)) = \prod_{j=1}^\ell    \Big (  \big( -2\sn(f_j) \hbar\big)^{\phi_+}\, c_\fg \, I_2\big(\downarrow,\sn(f_j)\big) \Big), \label{n2004}
\end{equation}
Dividing \eqref{n2007} by \eqref{n2004} and using Proposition \ref{n2005}, we have
\begin{align*} \hspace*{8pc}
\hWg(\hLMO(M))& = \prod_{j=1}^\ell |f_j|^{\phi_+}\, \tau^\fg(M)\\
&= |H_1(M,\Z)|\, \tau^\fg(M).
\end{align*}
This completes the proof the Theorem \ref{thm.rec} for the case when $M$ can be obtained by surgery along an algebraically split link.

Let us consider  the general case, when
$M$ is an arbitrary  rational homology 3-sphere.
It is known \cite{Ohtsuki_pi} that,
there exist some lens spaces $L(m_1,1)$, $\cdots$, $L(m_N,1)$
such that the connected sum
$M \# L(m_1,1) \# \cdots \# L(m_N,1)$
can be obtained from $S^3$ by surgery
along some algebraically split framed link.
Since the LMO invariant and the perturbative invariant are multiplicative
with respect to the connected sum,
it follows from the above case that
\begin{align*} \hspace*{3pc}
& \hWg \big( \hLMO (M) \big) \cdot
\prod_i \hWg \Big( \hLMO \big( L(m_i,1) \big) \Big)
\\
& = \
| H_1(M;\Z) |^{\phi_+} \,  \tau^\fg(M) \cdot
\prod_i \Big( \big| H_1 \big( L(m_i,1) ;\Z\big) \big|^{\phi_+} \,
\tau^\fg \big( L(m_i,1) \big) \Big) .
\end{align*}
In particular, since the lens space $L(m_i,1)$ can be obtained from $S^3$
by surgery along a framed knot,
it also follows from the above case that
$$
\hWg \Big( \hLMO \big( L(m_i,1) \big) \Big)
\, = \,
\big| H_1 \big( L(m_i,1) ;\Z\big) \big|^{\phi_+} \,
\tau^\fg \big( L(m_i,1) \big) .
$$
Further, since the leading coefficient of the LMO invariant is $1$,
the value of the above formula is non-zero.
Therefore, as the quotient of the above two formulas,
we obtain the required formula.
This completes the proof of Theorem \ref{thm.rec}
in the general case. \qed

%\subsection{$\Wg(\hLMO)$ as Gaussian integral on $\fg^*$}
\subsection{Some lemmas on weights of Jacobi diagrams}
\label{sec.lem_J}

In this section, we show some lemmas on Jacobi diagrams
which are used in the proof  of Proposition  \ref{n2005}.

\begin{lem}
\label{lem.W_U}
For a Jacobi diagram $D \in \cA(\ast)$ and a non-zero real number  $f$,
$$
\hWg \Big(
\Big\langle \exp \big( \frac{-1}{2f } \JdU{\ptl{x}}{\ptl{x}} \big),
D \Big\rangle \Big)
\, = \,
\eL^{(f)}_\fg \big( \hWg(D) \big)  \, .
$$
\end{lem}

\begin{proof}
The bracket can be presented in terms of differentials
as explained in \cite[Appendix]{Aint1}.
We verify this for the required formula concretely.

By expanding the exponential,
it is sufficient to show that
\begin{equation} \hspace*{8pc}
\label{eq.Ud_D}
\Wg \Big(
\Big\langle \big( \JdU{\ptl{x}}{\ptl{x}} \big)^d,
D \Big\rangle \Big)
\, = \,
\Dg^d \big( \Wg(D) \big) \big|_{X_i=0} \, .
\end{equation}
Since both sides are equal to $0$ unless $D$ has $2d$ legs,
we can assume that $D$ has $2d$ legs.

When $d=1$, (\ref{eq.Ud_D}) is shown by
$$\vspc{20}
\Wg \Big( \Big\langle
\JdU{\ptl{x}}{\ptl{x}} \!\! , \, D
\Big\rangle \Big)
\, = \,
\Wg \Big( 2 \,
\bp{45,10}{\put(0,0){\pc{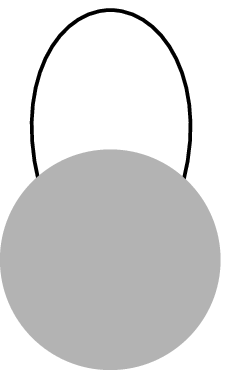}{0.5}}
\put(16,-10){$D$}} \Big)
\, = \,
2 \, B \big( \Wg(D) \big)
\, = \,
\Dg \big( \Wg(D) \big) ,
$$
where $B$ is the invariant form.
When $d=2$, putting\vspc{17}
$\Wg(D) = \sum_k Y_{1,k} Y_{2,k} Y_{3,k} Y_{4,k}$ for $Y_{i,j} \in \fg$,
(\ref{eq.Ud_D}) is shown by
\begin{align*} \hspace*{2pc} \vspc{29}
& \Wg \Big( \Big\langle
\big( \JdU{\ptl{x}}{\ptl{x}} \!\! \big)^2 , \, D
\Big\rangle \Big)
\, = \,
\sum_\tau
\Wg \Big(
\bp{61,10}{\put(0,0){\pc{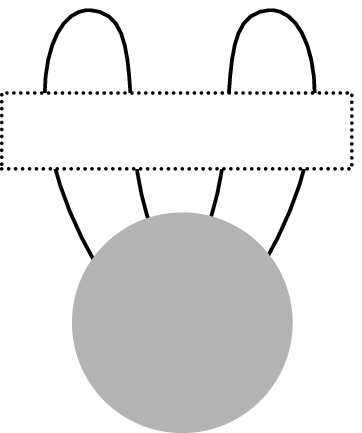}{0.5}}\put(27,15){$\tau$}\put(26,-14){$D$}} \Big)
\\ \vspc{20}
& \, = \,
\sum_{\tau,k}
B( Y_{\tau(1),k}, Y_{\tau(2),k} ) \, B( Y_{\tau(3),k}, Y_{\tau(4),k} )
\\
& \, = \,
\sum_{\tau,i,j,k}
\ptl{X_i}( Y_{\tau(1),k} ) \, \ptl{X_i}( Y_{\tau(2),k} ) \,
\ptl{X_j}( Y_{\tau(3),k} ) \, \ptl{X_j}( Y_{\tau(4),k} )
\, = \,
\Dg^2 \big( \Wg(D) \big) ,
\end{align*}
where the sum of $\tau$ runs over all permutations on $\{ 1,2,3,4 \}$.
For a general $d$, we can show (\ref{eq.Ud_D}) in the same way as above.
\end{proof}

\vskip 0.5pc

\vskip 0.5pc

\begin{lem}[\cite{Kuriya}]
\label{Wg_theta}
For the Jacobi diagram $\theta$ given in (\ref{eq.theta}), \
$\Wg(\theta) = 24 \, |\rho|^2$,
where $\rho$ is the half-sum of positive roots.
\end{lem}

\begin{proof}
It is shown from the definition of the weight system
(see, {\it e.g.}, \cite{Ohtsuki_book}) that
$$
\Wg \big( \pcp{27,10}{0,-2}{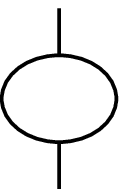}{0.5} \big) \, = \,
C_{\rm ad} \, \Wg \big( \pcp{20,10}{5,-2}{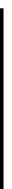}{0.5} \big)
\qquad \mbox{ and } \qquad
\Wg \big( \pcp{27,10}{-1,-2}{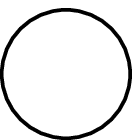}{0.5} \big) \, = \, {\rm dim} \, \fg ,
$$
where $C_{\rm ad}$ denotes the eigenvalue of the Casimir element
on the adjoint representation of $\fg$.
Hence, $\Wg(\theta) = C_{\rm ad} \, {\rm dim} \, \fg$.

It is known  that, $C_{\rm ad}= (\delta, \delta+2\rho)$, where
 $\delta$ is the highest weight of the adjoint representation, which is longest positive root. In our normalization of the inner product, $(\delta,\delta)=2d$ and
$(\delta,\rho) = d h^\vee - d  $,
where $h^\vee$ denotes the dual Coxeter number of $\fg$ and $d$ is the maximal absolute value of the off-diagonal entries of the Cartan matrix.
Therefore,
$\Wg(\theta) = 2 d h^\vee \, {\rm dim} \, \fg$.

Further, it is known \cite[47.11]{FdV} (adjusted to our normalization of the inner product) that
$2 d h^\vee \, {\rm dim} \, \fg = 24 |\rho|^2$.
Hence, we obtain the required formula.
\end{proof}

\section{The knot case}
\label{sec.knot}

The aim of this section is to prove Proposition \ref{n2003} for the case $\ell=1$. We call this the knot case, since Proposition \ref{n2003} with $\ell=1$ is enough to prove the main theorem for
the case when $M$ is obtained by surgery on a knot.
We show a geometric proof  in Section \ref{sec.geom_pf_knot}
and an algebraic proof  in Section \ref{sec.alg_pf_knot}.

\subsection{Geometric approach}
\label{sec.geom_pf_knot}

\begin{proof}[Geometric proof of Proposition \ref{n2003} in the knot case]
Since $\ell=1$, $g\in S(\fg)^\fg[[\hbar]]$. Without loss of generality, one can assume that $g\in S(\fg)^\fg$. We will write $f_1=f$.
Note that $g$ is a function on $\fg$; it's restriction on $\fh^*$ is denoted by $\cP(g)$. On the other hand, $\Ups_\fg(g) \in U(\fg)$ defines a function
on $\fh^*$,
see Section \ref{0001}.
{}From the commutativity of Diagram \ref{eq.cd}, we have that, as functions on $\fh^*$,
\begin{equation} \hspace*{8pc} \Ups_\fg(g) = \od \,\cP( g).
\label{n5000}
\end{equation}

 The left-hand side of \eqref{n2008} is $\eL^{(f)}_\fg(g)$, which, by Proposition \ref{n02}, can be expressed by an integral:
$$
\text{LHS of \eqref{n2008}} \, =\eL^{(f)}_\fg(g) = \frac{1}{(4\pi)^{\dim \fg/2}} \int_{\fg^*} e^{-|x|^2 /4}\, g (\frac{x}{\sqrt {-2f\hbar}} ) dx
$$
The integrand is invariant under the co-adjoint action.
Hence, according to Proposition \ref{prop.intg_inth} below,
one can reduce  the integral to an integral over
the Cartan subalgebra:
\begin{equation} \hspace*{4pc}
\text{LHS of \eqref{n2008}} \, \ = \
\frac{\tilde c_\fg}{(4\pi)^{\dim \fg/2}} \, \int_{\fh^*} \od^2(x)\,
e^{-|x|^2/4}
\cP(g) \Big(\frac{x}{\sqrt {-2f \hbar}}\Big)\, d x.
\label{n051}
\end{equation}
Here, $\tilde c_\fg$ is a non-zero constant depending on the Lie algebra $\fg$ only.

We turn to the right-hand side of \eqref{n2008}. Using \eqref{n5000} one has

$$ \text{RHS of \eqref{n2008}} \, =  \, c_\fg\,  (-2 f\hbar)^{\phi_+} \, \eL^{(f)}_\fh (\od^2 \cP(g)).$$
Again using Proposition \ref{n02} we have

\begin{equation} \text{RHS of \eqref{n2008}} \, =  c_\fg\,  (-2 f\hbar) \,  \frac{1}{(4\pi)^{\dim \fh/2}}\int_{\fh^*}
e^{-|x|^2/4}\, \od^2\Big(\frac{x}{\sqrt {-2f \hbar}}\Big)\,
g \Big(\frac{x}{\sqrt {-2f \hbar}}\Big)\, d x.
\label{n061}
\end{equation}

Because  $\od^2$ is a homogeneous polynomial of degree $2\phi_+$, one has
$$
\od^2(x) = (-2f\hbar)^{\phi_+}\, \od^2\Big(\frac{x}{\sqrt {-2f \hbar}}\Big).
$$
With $c_\fg= \frac{\tilde c_\fg}{(4\pi)^{\phi_+}}$, from \eqref{n051} and \eqref{n061} we see that
$$ \text{LHS of \eqref{n2008}} \, =\, \text{RHS of \eqref{n2008}}.$$
\end{proof}

\vskip 0.5pc

\subsubsection{Gaussian integral and $\eL^{(f)}_V$} Suppose $V$ is a Euclidean space and $f$ a non-zero number.
The following lemma says that the operator $\eL^{(f)}_V$ can be expressed by an integral.

\begin{lem} \label{n02} Suppose
$g \in S(V)(\!(h)\!)$, considered as a function on $V^*$ with values in $\BR(\!(h)\!)$.
Then
\begin{equation} \hspace*{8pc}
\eL^{(f)}_V(g)
\ = \
\frac{1}{(4\pi)^{\dim V/2}}\int_{V^*} e^{-  |x|^2/4}\, g\Big(
\frac{x}{\sqrt {-2f \hbar}}
\Big) \,  d x.
\label{n01}
\end{equation}
\end{lem}
\begin{rem} Here, $g\Big(
\frac{x}{\sqrt {-2f \hbar}}
\Big)$ is the  function on $V^*$ with values in $\C(\!( \hbar^{1/2})\!)$ defined as follows. If $g$ is of the form
$g = z^d$ where $z \in V$, then
$$ g\Big(
\frac{x}{\sqrt {-2f \hbar}}
\Big) := g(x)\, \left( (\sqrt {-2f \hbar})^{-d}\right) .$$
The square root in the right-hand side does not really appear, since if $d$ is odd, then  both sides of \eqref{n01} are 0.
\end{rem}
\begin{proof} We can assume that $g\in S(V)$.
Every polynomial is a sum of powers of  linear polynomials. Since both sides of \eqref{n01} depend linearly on $g$, we can assume that $g$ is a power of a linear polynomial.
By changing coordinates one can assume that $g = x_1^d$, where $x_1$ is the first of an orthonormal basis $x_1, \dots, x_n$ of $V$. The statement now reduces to the case when $V$ is one-dimensional, which follows from a simple Gaussian integral calculation,
see {\it e.g.} \cite[Lemma 2.11]{BGV}.
\end{proof}

\subsubsection{Reduction from $\fg^*$ to $\fh^*$}

\begin{prop}
\label{prop.intg_inth}
Suppose $g$ is a  $G$-invariant  function on $\fg^*$. Then
$$
\int_{\fg^*} g \, d x \ = \ \tilde c_\fg \int_{\fh^*} \od ^2 \, \cP(g) \, d x
$$
provided that both side converges absolutely.
Here, $\tilde c_\fg$ is a non-zero constant depending only on $\fg$.
\end{prop}

\begin{proof}
It is clear that if such $\tilde c_\fg$ exists, then it is non-zero,
since there are $G$-invariant functions $g$, {\it e.g.}
$g(x) = \exp (-|x|^2)$,  for which the left-hand side is non-zero.

The co-adjoint action of $G$ on $\fg^*$ is well-studied in the literature. A point $x \in \fg^*$ is {\em regular} if its orbit $G\cdot x$ is a submanifold of  dimension $\dim \fg -\dim \fh= 2\phi_+$, the maximal possible dimension.
It is known that the set of non regular points has measure 0. Every orbit has non-empty intersection with $\fh^*$, and
if $x$ is regular, then $G\cdot x \cap \fh^*$ has exactly $|W|$ points. Since the function $g$ is constant on each orbit, we have
$$
\int_{\fg^*}  g(x) d x \ = \
\frac{1}{|W|} \int_{\fh^*} \vol(G\cdot x) \, \cP(g)(x)  d x .
\label{1}
$$

The volume function is also well-known;
it can be calculated, for example, from \cite[Chapter 7]{BGV}:
\begin{equation} \hspace*{12pc}
\vol(G \cdot x) \, = \, \tilde c_\fg'\,  \od^2(x)
\label{030}
\end{equation}
where $\tilde c_\fg'$ is a constant.
{}From \eqref{030} we can deduce the proposition, with $\tilde c_\fg =\tilde  c'_\fg/|W|$.

Here is a simple proof\footnote{
The authors thank A. Kirillov Jr. for supplying them the proof.}
of \eqref{030}.
We will identify $\fg$ with $\fg^*$ via the invariant inner product.
Let $H$ be the maximal abelian subgroup of $G$ whose Lie algebra is $\fh$. The space $G/H$ is a homogeneous $G$-space. The tangent space of $G/H$ at $H$ can be identified with
$\fh^\perp$, with inner product induced from the invariant; from this we define a Riemannian metric on $G/H$.
When $x\in \fh$ is a regular, its stationary group is isomorphic to the torus $H$.
The map $\varphi: G/H \to G\cdot \fh$, defined by $g \to g\cdot x$ with $g\in G$, is a diffeomorphism.
 The tangent space of $G\cdot x$ at $x$ can also be identified with
the same $\mathfrak \fh^\perp$ with the same inner product.
It is easy to see that $\varphi$ at $H$ has derivative $d\varphi_H=-\ad(x): \fh^\perp \to \fh^\perp$. Let us calculate the
determinant of $d\varphi$. Because $G/H$ is $G$-homogeneous and $\varphi$ is $G$-equivariant, $|\det(d \varphi)|$ is constant on $G/H$, hence $ |\det(d \varphi)| = |\det(\ad (x)|$.
To calculate  $|\det(\ad (x)|$, it's easier to use the complexification of the adjoint representation, since $\ad(x)$ is diagonal in the complexified representation. The complexified
 $\fh^\perp_\C$ has the standard Chevalley  basis $E_\al, F_\al, \al\in \Phi_+$ such that  $\ad(x) E_\al= i (x,\al) E_\al$ and $\ad(t) F_\al = - i(x,\al) F_\al$. It follows that
$|d \varphi|= \prod_{\al \in \Phi_+} |(x,\al)|^2$. Hence
$$
\vol(G\cdot x) \ = \ \vol(G/H) \, \prod_{\al \in \Phi_+} |(x,\al)|^2
\ = \ \tilde c'_\fg \, \od^2(x),
$$
where $\tilde c'_\fg= \vol(G/H) \, \prod_{\al \in \Phi_+} |(\rho,\al)|^2$.
\end{proof}

\subsection{Algebraic approach}
\label{sec.alg_pf_knot}

In this section, we show an algebraic proof of Proposition \ref{n2003}
in the knot case, {\it i.e.} the case $\ell=1$.
We also verify some formulas of the proof
in the $\fsl2$ case and in the $\fsl3$ case
in Sections \ref{sec.sl2} and \ref{sec.sl3} respectively.

\begin{proof}[Algebraic proof of Proposition \ref{n2003} in the knot case]
Again we can assume that $g \in S(\fg)^\fg$.

By definition, the left-hand side of \eqref{n2008} is
$$
\eL^{(f)}_\fg(g) \, = \,
 \exp \big( - \frac{1}{2f \hbar} \, \Dg \big)
\big( g \big) \Big|_{x = 0} \, .
$$
By expanding the exponential,
\begin{equation}
\label{eq.sum_Dd_g} \hspace*{7pc}
\text{LHS of \eqref{n2008}} \,  \, = \,
 \sum_{d \ge 0} \big( - \frac{1}{2 f \hbar} \big)^d \,
\frac{1}{d \, !} \, \Dg^d (g_d),
\end{equation}
where $g_d$ is the degree $2d$ part of $g$.

\vskip 0.2pc

Let us turn to the right-hand side of \eqref{n2008}. Recall that $\od$ has degree $\phi_+$. By \eqref{n5000}

\begin{align} \hspace*{1pc} \notag
\text{LHS of \eqref{n2008}} \,
& \, = \, c_\fg (-2f\hbar)^{\phi_+} \, \eL^{(f)}_\fh (\od^2 \cP(g)) \\
&\, = \, c_\fg (-2f\hbar)^{\phi_+} \,
\exp \big( - \frac{1}{2f \hbar} \Dh \big)
\big(
\od^2 \,
\cP(g)  \big)
\Big|_{x=0} \notag
\\
& \, = \, c_\fg (-2f\hbar)^{\phi_+} \,
 \sum_{d \ge 0}
\big( - \frac{1}{2 f \hbar} \big)^{d+\phi_+}
\frac{1}{(d+ \phi_+)!}
\Dh^{d+ \phi_+} \big( \od^2 \, \cP(g_d) \big) \notag\\
& \, = \, c_\fg  \,
 \sum_{d \ge 0}
\big( - \frac{1}{2 f \hbar} \big)^{d}
\frac{1}{(d+ \phi_+)!}
\Dh^{d+ \phi_+} \big( \od^2 \, \cP(g_d) \big).
\label{eq.sum_Dh_g}
\end{align}
%\label{eq.b}

\vskip 0.2pc

Comparing (\ref{eq.sum_Dd_g}) and (\ref{eq.sum_Dh_g})
by using Proposition \ref{prop.DhD} below, we have immediately
$$  \text{LHS of \eqref{n2008}} \, = \text{RHS of \eqref{n2008}} \, .  $$
This completes the algebraic proof of Proposition \ref{n2003} in the knot case.
\end{proof}

\vskip 0.5pc

\begin{prop}
\label{prop.DhD}
For any homogeneous polynomial
$g \in S(\fg)^\fg$ of degree $2d$,
$$
\frac{c_\fg}{d \, !} \, \Dg^d (g)
\ = \ \frac{1}{(d+ \phi_+)!} \, \Dh^{d+ \phi_+} \big( \od^2 \, \cP(g )\big),
$$
where $c_\fg$ is a non-zero constant depending on $\fg$ only.
\end{prop}

\begin{proof} Since the right hand side is not identically 0, if such a $c_\fg$ exists, then it is non-zero.
We show that the identity of the proposition holds true if we take $c_\fg = \Dh^{\phi_+} \big( \od^2 \, \big)/(\phi_+)!$.

Since $\Dg^d(g)$ is a scalar, we have that
$$
\od \, \Dg^d(g)
= \od \, \cP \big( \Dg^d(g) \big)
= \Dh^d \big( \od \, g \big),
$$
where we obtain the second equality
by applying Proposition \ref{prop.HCrf} below repeatedly.
Hence, substituting the above formula,
$$
\Dh^{\phi_+} \big( \od^2 \, \Dg^d(g) \big)
\, = \, \Dh^{\phi_+} \Big( \od \, \Dh^d \big( \od \, \cP(g) \big) \Big).
$$
Further, since the left-hand side is presented by
$\Dh^{\phi_+} \big( \od^2 \, \Dg^d(g) \big)
= \Dh^{\phi_+} (\od^2) \, \Dg^d(g) = $ \mbox{$c_\fg \cdot \phi_+ ! \, \Dg^d(g)$},
the required formula is reduced to
$$
\Dh^{d+\phi_+} \big( \od^2 \, \cP(g) \big)
\, = \, \binom{d+\phi_+}{\phi_+} \,
\Dh^{\phi_+} \Big( \od \, \Dh^d \big( \od \, \cP(g) \big) \Big).
$$
It is sufficient to show this formula.

By putting $g' = \od \, \cP(g)$, the above formula is rewritten,
$$
\Dh^{d+\phi_+} \big( \od \, g' \big) \, = \,
\binom{d+\phi_+}{\phi_+} \, \Dh^{\phi_+} \big( \od \, \Dh^d ( g' ) \big).
$$
As for $\Dh^{d+\phi_+}$ in the left-hand side,
since $\Dh(\od) = 0$ by Lemma \ref{lem.Db=0} below,
$\phi_+$ copies of $\Dh$ in $\Dh^{d+\phi_+}$ act on $\od$.
The number of choices of these $\phi_+$ copies is
the binomial coefficient in the right-hand side.
Further, these $\phi_+$ copies with $\od$ can be replaced by
a differential operator with scalar coefficients,
and this differential operator commutes $\Dh$.
Hence, we obtain the above formula.
\end{proof}

\vskip 0.5pc

\begin{lem}
\label{lem.Db=0}
One has
$\Dh(\od) = 0$.
\end{lem}

\begin{proof} Note that $\od$ is a $W$-anti-symmetric polynomial. Actually, there is no $W$-anti-symmetric polynomial of degree lower than that of $\od$. This follows from the well-known fact that
the ring of  $W$-anti-symmetric polynomials on $\fh^*$ is the free module over the ring of $W$-symmetric polynomials spanned by $\od$.
Since $\Dh$ is $W$-invariant, $\Dh(\od)$ is a $W$-anti-symmetric polynomial of lower degree, hence it is 0.
\end{proof}

\vskip 0.5pc

\vskip 0.5pc

The following proposition is a reformulation of
Harish-Chandra's restriction formula;
see \cite[Proposition II.3.14]{Helgason}, \cite[Theorem 2.1.8]{HT}.

\begin{prop}
\label{prop.HCrf}
For any $g \in S(\fg)^\fg$
$$
\od \, \cP \big( \Dg(g) \big)
\, = \, \Dh \big( \od \, \cP(g) \big) .
$$
\end{prop}

\begin{proof}
By using the invariant form, we identify $\fg$ and $\fg^*$,
and identify $\fh$ and $\fh^*$.
Further, $\od \in P(\fh^*)$ is identified with $\od^* \in P(\fh)$ given by
$$
\od^*(X) \, = \, \prod_{\al \in \Phi_+} \frac{\al(X)}{(\rho,\al)}
$$
for $X \in \fh$.
Then, the required formula is identified with
Harish-Chandra's restriction formula
(see \cite[Proposition II.3.14]{Helgason}, \cite[Theorem 2.1.8]{HT}).
\end{proof}

\subsubsection{The $\fsl2$ case}
\label{sec.sl2}

Although we have proved Propositions \ref{prop.DhD},  we will write down explicitly the identity of this proposition in the case of $sl_2$ (and $sl_3$ in the next section), and verify
Propositions \ref{prop.DhD} and \ref{prop.HCrf} by direct calculation. The reader will see that the identity is quite non-trivial.

We recall that $\gC$ is spanned by
$$
H = {\normalsize \begin{pmatrix} 1 & 0 \\ 0 & -1 \end{pmatrix}},
\qquad
E = {\normalsize \begin{pmatrix} 0 & 1 \\ 0 & 0 \end{pmatrix}},
\qquad
F = {\normalsize \begin{pmatrix} 0 & 0 \\ 1 & 0 \end{pmatrix}}.
$$
According to the convention, we give the invariant form
by $(X,Y) = {\rm Tr}( X Y )$.
We regard $H$, $E$, $F$ as variables in the following of this section.
The Cartan subalgebra $\fh$ is spanned by $H$.
We choose the fundamental weight
$\La \in \fh^*$ such that $\La(H)=1$.
Since $|H|^2 =2$, \ $|\La|^2 = \frac12$.
Further, since $\Phi_+ = \{ 2 \La \}$, \ $\rho = \La$.
We put $\la = n \La$.
We denote by $V_n$ the $n$-dimensional irreducible representation,
which is the irreducible representation whose highest weight is $\la - \rho$.
Since the Laplacian $\Dg$ is characterized by
$\frac12 \Dg(X Y) = (X,Y)$,
it is presented by
$\Dg \, = \, 2(\ptl{H}^2 + \ptl{E} \ptl{F})$,
which acts on $S(\gC) = \C[H,E,F]$.
Similarly, the Laplacian $\Dh$ is presented by
$\Dh = 2 \, \ptl{n}^2$, which acts on $P_\C(\fh^*)^W = \C[n^2]$.
It is known that
$S(\gC)^\fg = \C[C]$,
where $C$ is the Casimir element given by
$C = \frac12 H^2 + 2 E F$.

When $\fg = \fsl2$, Proposition \ref{prop.DhD} is rewritten,
$$
\frac{4}{d \, !} \, \Dg^d (C^d)
\ = \ \frac{1}{(d+1)!} \, \Dh^{d+1} \big( n^2 \, \cP(C^d) \big) ,
$$
where $\cP : \C[C] \lra \C[n^2]$ is given by $\cP(C) = n^2/2$.
Since
$\Dg^d(C^d) = \frac{(2d+1)!}{2^d}$
and
$\Dh^{d+1} \big( n^2 \, \cP(C^d) \big) = \frac{(2d+2)!}{2^d}$,
we can verify Proposition \ref{prop.DhD} for $\fg = \fsl2$,
by calculating both sides of the above formula concretely.

\vskip 0.2pc

When $\fg = \fsl2$, Proposition \ref{prop.HCrf} is rewritten,
$$
n \, \cP \big( \Dg(g) \big)
= \Dh \big( n \, \cP(g) \big)
$$
for $g \in S(\gC)^\fg = \C[C]$.
We verify this formula, as follows.
It is sufficient to show the formula when $g = C^d$.
Then,
$\Dg(C^d) = d (2d+1) C^{d-1}$ and
$n \, \cP \big( \Dg(C^d) \big) = \frac{d(2d+1)}{2^{d-1}} n^{2d-1}$.
On the other hand,
$\Dh \big( n \, \cP(C^d) \big) = \Dh \big( \frac{n^{2d+1}}{2^d} \big)
= \frac{d(2d+1)}{2^{d-1}} n^{2d-1}$.
Hence, Proposition \ref{prop.HCrf} was verified for $\fg = \fsl2$.

\subsubsection{The $\fsl3$ case}
\label{sec.sl3}

In this section,
we give some calculation to verify
Propositions \ref{prop.DhD} and \ref{prop.HCrf} concretely,
when $\fg = \fsl3$.

We recall that $\gC$ is spanned by
\begin{align*} \hspace*{1pc}
& H_1 = {\normalsize \begin{pmatrix}
1\!&0&0 \\ 0&\!-1\! &0 \\ 0&0&\!0\end{pmatrix}},
\ \
H_2 = {\normalsize \begin{pmatrix}
0 &0&0 \\ 0& 1\! &0 \\ 0&0& \!\!-1 \end{pmatrix}},
\ \
E_1 = {\normalsize \begin{pmatrix}
0 & 1 &0 \\ 0& 0 &0 \\ 0&0& 0 \end{pmatrix}},
\ \
E_2 = {\normalsize \begin{pmatrix}
0 &0&0 \\ 0& 0 &1 \\ 0&0& 0 \end{pmatrix}},
\\
& E_3 = {\normalsize \begin{pmatrix}
0 &0&0 \\ 0& 0 &0 \\ 1&0& 0 \end{pmatrix}},
\ \
F_1 = {\normalsize \begin{pmatrix}
0 &0&0 \\ 1& 0 &0 \\ 0&0& 0 \end{pmatrix}},
\ \
F_2 = {\normalsize \begin{pmatrix}
0 &0&0 \\ 0& 0 &0 \\ 0&1& 0 \end{pmatrix}},
\ \
F_3 = {\normalsize \begin{pmatrix}
0 &0& 1 \\ 0& 0 &0 \\ 0&0& 0 \end{pmatrix}}
\end{align*}
According to the convention, we give the invariant form
by $(X,Y) = {\rm Tr}( X Y )$.
We regard $H_1, H_2, \cdots$  as variables in the following of this section.
The Cartan subalgebra $\fh$ is spanned by $H_1$ and $H_2$.
We choose the fundamental weights
$\La_1,\La_2 \in \fh^*$ such that
$\La_i(H_i) = 1$ for $i=1,2$, and $\La_i(H_j) = 0$ if $i \ne j$.
Further,
since $\Phi_+ = \{ 2 \La_1 -\La_2, \, \La_1+\La_2, \, 2 \La_2-\La_1 \}$,
\ $\rho = \La_1 + \La_2$.
We put $\la = n \La_1 + m \La_2$.
We denote by $V_{n,m}$
the irreducible representation of $\fsl3$
whose highest weight is $\la - \rho$,
which is the irreducible representation presented by
the following Young diagram.
$$
\begin{picture}(110,25)
\put(0,5){\line(1,0){55}}\put(0,15){\line(1,0){110}}\put(0,25){\line(1,0){110}}
\put(0,5){\line(0,1){20}}\put(10,5){\line(0,1){20}}\put(20,5){\line(0,1){20}}
\put(45,5){\line(0,1){20}}\put(55,5){\line(0,1){20}}\put(65,15){\line(0,1){10}}
\put(75,15){\line(0,1){10}}\put(100,15){\line(0,1){10}}
\put(110,15){\line(0,1){10}}
\put(25,7){$\cdots$}\put(25,17){$\cdots$}\put(80,17){$\cdots$}
\put(57,13){\normalsize $\underbrace{\hspace*{4.3pc}}_{n-1}$}
\put(1,3){\normalsize $\underbrace{\hspace*{4.3pc}}_{m-1}$}
\end{picture}
$$
Since the Laplacian $\Dg$ is characterized by
$\frac12 \Dg(X Y) = (X,Y)$,
it is presented by
$$
\Dg \, = \, 2 \big( \ptl{H_1}^2 + \ptl{H_2}^2 - \ptl{H_1} \ptl{H_2}
+ \ptl{E_1} \ptl{F_1} + \ptl{E_2} \ptl{F_2} + \ptl{E_3} \ptl{F_3} \big) ,
$$
which acts on
$S(\gC) \, = \, \C \big[ H_1, H_2, E_1, E_2, E_3, F_1, F_2, F_3 \big]$.
Similarly, the Laplacian $\Dh$ is presented by
$$
\Dh \, = \, 2 \big( \ptl{n}^2 + \ptl{m}^2 - \ptl{n} \ptl{m} \big) ,
$$
which acts on $P_\C(\fh^*)^W = \C[n,m]^W$.
The map
$$
\cP : \, \C[H_1,H_2,E_1,E_2,E_3,F_1,F_2,F_3]^\fg \lra \C[n,m]^W
$$
is given by
$\cP(H_1)=n$, \, $\cP(H_2)=m$, \, $\cP(E_i)=\cP(F_i)=0$.
Further, $\od(\la) = n m (n+m)/2$.

When $\fg = \fsl3$,
Proposition \ref{prop.DhD} is rewritten,
\begin{equation} \hspace*{8.5pc}
\label{eq.DhD_sl3}
\frac{24}{d \, !} \, \Dg^d (g)
\ = \ \frac{1}{(d+3)!} \, \Dh^{d+3} \big( \od^2 \, \cP(g) \big)
\end{equation}
for any homogeneous polynomial $g \in S(\gC)^\fg$ of degree $2d$.
Further,
Proposition \ref{prop.HCrf} is rewritten,
\begin{equation} \hspace*{9.5pc}
\label{eq.HCrf_sl3}
\od \, \cP \big( \Dg(g) \big)
\, = \, \Dh \big( \od \, \cP(g) \big)
\end{equation}
for any $g \in S(\gC)^\fg$.
It is known that $S(\gC)^\fg = \C[C, C_3]$, where
\begin{align*} \hspace{2pc}
& C \, = \, \frac13 (H_1^2 +H_2^2 +H_1 H_2) +E_1 F_1 +E_2 F_2 +E_3 F_3
\\
& C_3 \, = \, - \frac19 (H_1-H_2)(H_2-H_3)(H_3-H_1)
+ 3 E_1 E_2 E_3 + 3 F_1 F_2 F_3
\\*
& \qquad\quad + E_1 E_1 (H_2-H_3) + E_2 E_2 (H_3-H_1) + E_3 F_3 (H_1-H_2),
\end{align*}
putting $H_3 = - H_1 - H_2$.
By computer calculation,
we can verify (\ref{eq.DhD_sl3}) and (\ref{eq.HCrf_sl3})
for concrete $g \in S(\gC)^\fg$ of small degrees
putting $g$ to be polynomials in $C$ and $C_3$.

\section{The link case}
\label{sec.link}

In Section \ref{sec.knot},
we gave proofs of Proposition \ref{n2003}, and hence Theorem \ref{thm.rec}, in the knot case.
Here we give a proof of  Proposition \ref{n2003} in the general case in Section \ref{sec.link_alg}.
In Section \ref{sec.link_fti},
we also show that,
without Proposition \ref{n2003} for the case $\ell >1$, one can still
prove the main theorem
 using general results on finite type invariants.

\subsection{The link case by direct calculation}
\label{sec.link_alg}

\begin{proof}[Proof of Proposition \ref{n2003} in the link case]
 The left-hand side of \eqref{n2008} is
$$
\text{LHS of \eqref{n2008}} \, = \, \eL^{(\ff)}_\fg(g) = \left (\bigotimes_{j=1}^\ell \eL^{(f_j)}_\fg \right) (g).
$$
Note that $\eL^{(f)}_\fh$ acts on $P(\fh^*)$.
We define a modification of $\eL^{(f)}_\fh$, which acts on the bigger space $ P(\fg^*)=S(\fg)$, as follows:
\begin{equation} \hspace*{8pc}
\teL^{(f)}_\fh (g) :=
\left (-2 f_j \hbar)^{\phi_+}\, c_\fg\right)
\eL^{(f)}_\fh( \od \, \Ups_\fg (g)).
  \label{n3000}
\end{equation}
Then the right-hand side of \eqref{n2008} can be rewritten as
$$
\text{RHS of \eqref{n2008}} \, = \,  \left (  \bigotimes_{j=1}^\ell \teL^{(f_j)}_\fh \right )(g).
$$
Proposition \ref{n2003} becomes the statement that
for any $g \in (S(\fg)^{\otimes \ell })^\fg$,
$$
\left (\bigotimes_{j=1}^\ell \eL^{(f_j)}_\fg \right) (g)\,  = \, \left (  \bigotimes_{j=1}^\ell \teL^{(f_j)}_\fh \right )(g),
$$
which is the case $m=\ell$ of the following identity,
\begin{equation} \hspace*{7pc}
\label{eq.hcDg} \Big( \bigotimes_{1 \le j \le m} \!\! \eL^{(f_j)}_\fg
\otimes \bigotimes_{m < j \le \ell}  \!\! \teL^{(f_j)}_\fh \Big) (g)
\ = \
\,   \bigotimes_{j=1}^\ell \teL^{(f_j)}_\fh(g) .
\end{equation}

We will prove \eqref{eq.hcDg} by induction on $m$. The case $m=0$ is a tautology. Note also that when $\ell=1$, the identity holds since we proved it
in Section \ref{sec.knot}.
We put
$$
g' \, = \,
\Big( \bigotimes_{1 \le j < m} \!\! \eL^{(f_j)}_\fg
\otimes \bigotimes_{m < j \le \ell}  \!\! \teL^{(f_j)}_\fh \Big)(g)
\ \in \ S(\fg).
$$
Then equality \eqref{eq.hcDg} becomes
\begin{equation} \hspace*{12pc}
 \eL^{(f_m)}_\fg (g') = \teL^{(f_m)}_\fh(g')\ .
 \label{n3003}
\end{equation}
Since $\eL^{(f_j)}_\fg$ and $\teL^{(f_j)}_\fh $ are intertwiners by
Lemmas \ref{lem.cDf} and \ref{lem.hcDf} below,
$g' \in S(\fg)^\fg$.
Hence, \eqref{n3003} follows from the case $\ell=1$,
completing the induction.
\end{proof}

\vskip 0.5pc

\vskip 0.5pc

\begin{lem}
\label{lem.cDf} The map
$\eL^{(f)}_\fg : S(\fg)[[\hbar]] \to \R(\!(\hbar)\!)$ is an intertwiner
with respect to the action of $\fg$.
\end{lem}

\begin{proof}
By definition,
$\eL^{(f)}_\fg$ takes a monomial of odd degree in $S(\fg)$ to $0$. It is enough to consider the case $g =Y_1 Y_2 \cdots Y_{2d}$ where each $Y_j$ is a linear form.
Then $\eL^{(f)}_\fg$ takes
$Y_1 Y_2 \cdots Y_{2d}$ to a constant multiple of
$$
\sum_\tau B(Y_{\tau(1)},Y_{\tau(2)}) \cdots B(Y_{\tau(2d-1)},Y_{\tau(2d)}),
$$
where the sum runs over all permutations on $\{1,2, \cdots, 2d\}$ and $B$ is the invariant inner product.
Since the invariant form $B$ is an intertwiner,
$\eL^{(f)}_\fg$ is also an intertwiner.

Another proof is to use Proposition \ref{n02} to present $\eL^{(f)}_\fg$ by an integral:
$$ \eL^{(f)}_\fg (g) = \frac{1}{(4\pi)^{\dim \fg/2}}\int_{\fg^*}
e^{- |x|^2/2}\,  g \left( \frac{x}{\sqrt{-2 f \hbar}}  \right) dx\ .$$
Since $|x|^2$ and $dx$ are $G$-invariant, the right-hand side is $G$-invariant.
\end{proof}

\vskip 0.5pc

\begin{lem}
\label{lem.hcDf}
The map
$\teL^{(f)}_\fh : S(\fg)[[\hbar]] \to \R(\!(\hbar)\!)$ is an intertwiner
with respect to the action of $\fg$.
\end{lem}

\begin{proof}
Since the $\fg$ acts trivially  on $\R$,
it is sufficient to show that
$$
\teL^{(f)}_\fh \big(  {\rm ad}_X  (g) \big) = 0
$$
for $X \in \fg$ and $g \in S(\fg)$. Using the definition of $\teL^{(f)}_\fh$ in \eqref{n3000}, this is equivalent to
$$ \eL^{(f)}_\fh \Big( \od \, \Ups_\fg\left( {\rm ad}_X (g)  \right)\Big) =0.$$
It is enough to show that $\Ups_\fg\left( {\rm ad}_X (g)  \right) =0$ as a function on $\fh^*$. Evaluating $\Ups_\fg\left( {\rm ad}_X (g)  \right)$ on $\lambda \in \fh^*$ such that
$\lambda -\rho$ is a dominating weight, one has
\begin{align*}
\Ups_\fg\left( {\rm ad}_X (g)  \right) (\lambda) & = \Tr_{V_{\la-\rho}} \Ups_\fg \left( {\rm ad}_X (g)  \right) \quad \text{by definition } \\
 &= \Tr_{V_{\la-\rho}}   {\rm ad}_X  \left( \Ups_\fg(g)  \right) \quad \text{since $\Ups_\fg$ is an intertwiner}\\
 &= \Tr_{V_{\la-\rho}}\left(     X \, \Ups_\fg(g)  - \Ups_\fg(g)\, X \right) \quad \text{by definition of ${\rm ad}_X$ on $U(\fg)$} \\
 &=0 .
 \end{align*}
\end{proof}

\subsection{The link case through the knot case}
\label{sec.link_fti}

Here we give another approach to the link case using general results on finite type invariants. We will prove that if two multiplicative finite type invariants of rational homology 3-spheres coincide on the set of rational homology 3-spheres obtained by surgery on knots, then they are equal.

Let $\cH_1$ be the set of all integral homology 3-spheres which can be obtained by surgery along knots with framing $\pm1$, and $\cH_1^\oplus$ the set of all finite connected sums of elements in $\cH_1$.

\subsubsection{Finite type invariants of rational homology 3-spheres} We summarize here some basic facts about finite type invariants of {\em rational} homology 3-spheres (Ohtsuki, Goussarov-Habiro, for details see \cite{GGP,Ha2}).

\begin{figure}[htb]
$$ \psdraw{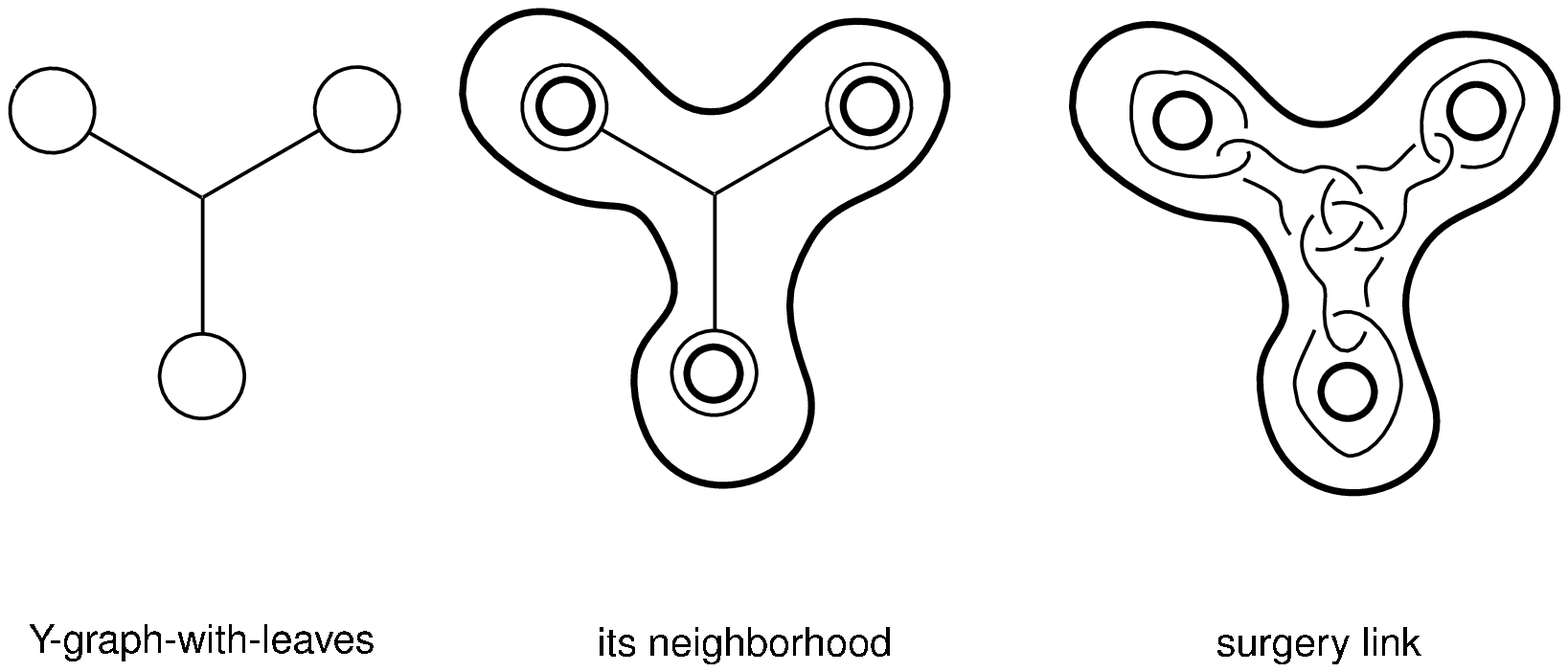}{4in} $$
\caption{\label{13}}
\end{figure}

Consider the {\em
standard $Y$-graph} in $\R^3$, see Figure \ref{13}.
 A $Y$-graph $C$ in $M$ is the image of an
embedding of a small neighborhood of the standard
$Y$-graph into $M$. Let $L$ be the six-component link
in a small neighborhood of the standard $Y$-graph as
shown in Figure \ref{13}, each component having framing 0. The surgery of
$M$ along the image of the six-component link is called a
$Y$-surgery along $C$, denoted by $M_C$.

Matveev \cite{Matveev} proved that  $M$ and $M'$ are related by a finite sequence of $Y$-surgeries if
and only if there is an isomorphism from $H_1(M,\Z)$ onto $H_1(M',\Z)$ preserving
the linking form on the torsion group. For a 3-manifold $M$ let $\cM(M)$ be the free $R$-module with basis
all 3-manifolds which have the same $H_1$ and linking form as $M$. Here $R$ is a commutative ring with unit. For example, $\cM(S^3)$ is the free $R$-module spanned by all integral homology 3-spheres.
We will always assume that $2$ is invertible in $R$. Actually, for the application in this paper, it's enough to consider the case when $R$ is a field of characteristic 0.

Let $E$ be a finite collection of disjoint
$Y$-graphs in a 3-manifold $N$. Define
$$[N,E] = \sum_{E'\subset E} (-1)^{|E'|} N_{E'}.$$

Define $\cF_n \cM(M)$
as $R$-submodule of $\cM(M)$ spanned by all $[N,E]$ such that $N$ is in
$\cM(M)$ and $|E|=n$. Any invariant $I$ of 3-manifolds in $\cM(M)$ with values in an $R$-module $A$ can be extended linearly to an $R$-linear function $I: \cM(M) \to A$.
Such an invariant $I$ is a {\em finite type invariant of order} $\le n$ if $I|_{ \cF_{n+1}} =0$. Matveev's result shows that an invariant of degree 0 is a constant
invariant in each class $\cM(M)$.

\begin{figure}[htb]
$$ \psdraw{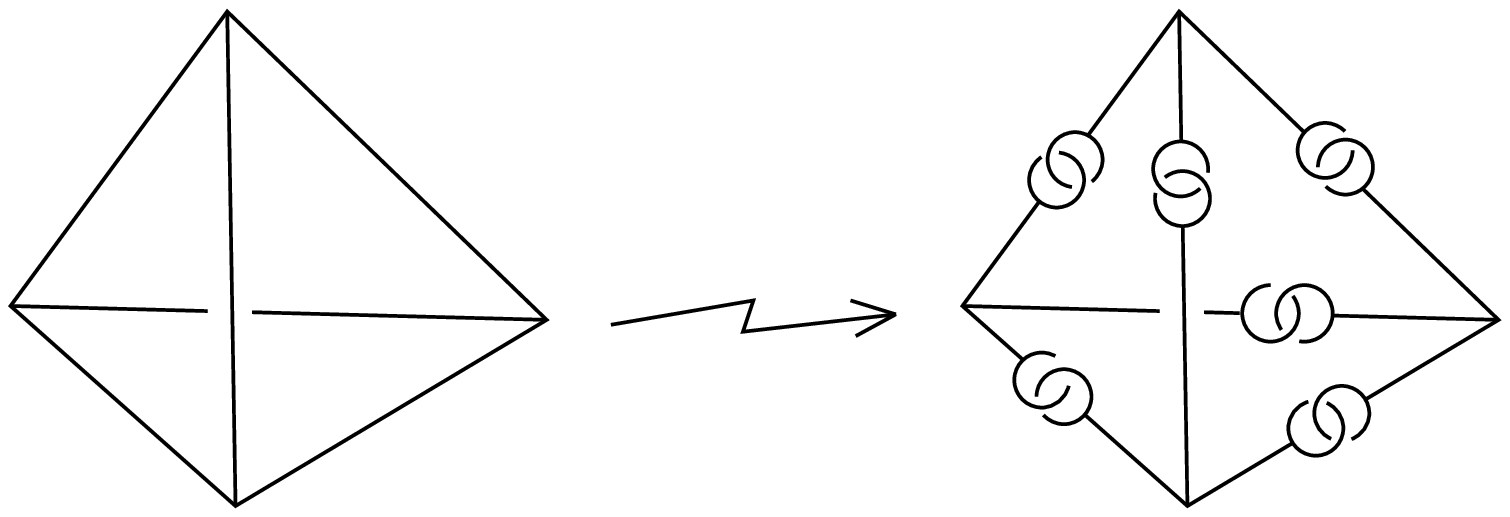}{3.6in} $$
\caption{\label{15}}
\end{figure}

Goussarov and Habiro showed that $\cF_{2n-1} = \cF_{2n}$. There is a {\em surjective} map
$$W: \Gr_n \cA(\emptyset)\to
\cF_{2n}\cM(M)/\cF_{2n+1}\cM(M),$$
 known as the universal weight map,  defined as follows.
Suppose $D\in \Gr_n \cA(\emptyset)$ is a Jacobi graph of degree
$n$.  Embed $D$ into $S^3$ arbitrarily.
Then from the image of $D$ construct a set $E$ of
$Y$-graphs as in  Figure \ref{15}.

By definition, $[M\#S^3,E]\in \cF_{2n}\cM(M)$. A priori, $[M\#S^3,E]$ depends on the way $D$ is embedded in $S^3$. However,

$$W(D) := [M\#S^3,E] \pmod {\cF_{2n+1}\cM(M)}$$
depends only on $D$ as an element in $\cA(\emptyset)$. Moreover, the map $W: \Gr_n \cA(\emptyset) \to
\cF_{2n}\cM(M)/\cF_{2n+1}\cM(M)$, known as the universal weight, is  {\em surjective}.

\label{019}
\begin{lem} Suppose $D$ is connected. Then $S^3_E$ can be obtained by surgery on $S^3$ along a knot with framing $\pm1$, $S^3_E \in \cH_1$.
\label{021}
\end{lem}

\begin{proof}
Choose a sublink $E'$ of $E$ consisting of all components of $E$ except for one component $K$, and do surgery along this sublink. Using repeatedly the move which removes a zero-framing trivial knot together with another knot piercing the trivial knot, it is easy to see that the resulting manifold is still $S^3$. Let $K'$ be the image of $K$ is the resulting $S^3$. Now one has $S^3_E= S^3_{K'}$, an integral homology 3-sphere. The framing of $K'$ must be $\pm 1$ because the resulting is
an integral homology 3-sphere.
\end{proof}

If $I$ is a finite type invariant of degree $\le 2n$, then its $n$-th weight is defined as
the composition
$$ w_I^{(n)} = I \circ W : \Gr_n \cA(\emptyset) \to V.$$
It is clear that if $w_I^{(n)}=0$, then $I$ has degree $\le 2n-2$.
\subsubsection{Multiplicative finite type invariants and surgery on knots} The following result shows that finite invariants are determined by their values on a smaller subset of the set of
all applicable 3-manifolds. Besides application to the proof of the LMO conjecture, the result is also interesting by itself.

\begin{thm} \label{022}
\mbox{}
\\
{\rm a)} Suppose $I$ is a finite type invariant of
{\em integral} homology 3-spheres with values in an $R$-module $A$
such that $I(M)=1$ for every $M\in \cH_1^\oplus$.
Then $I(M)=1$ for every integral homology 3-sphere.
\\
{\rm b)} Suppose $I$ is a {\em multiplicative}  finite type  invariant
of {\em rational} homology 3-spheres with values in an $R$-algebra $A$. If
 $I(M)=1$ for every $M\in \cH_1$ and every lens space $M= L(p,1)$,
then $I(M)=1$ for every rational homology 3-sphere.
In particular, if $I(M)=1$ for any rational homology 3-sphere obtained by surgery on knots, then $I(M)=1$ for any rational homology 3-sphere.
\end{thm}

\begin{proof}

a) Suppose $I$ has degree $\le 2n$. Let $D$ be a Jacobi diagram of degree $n$.  Suppose $D= \prod_{j=1}^s D_j$. Let $E_j$ be the $Y$-graphs corresponding to  $D_j$ as constructed in Subsection \ref{019}, and $E= \sqcup_{j=1}^s E_j$. Since each of $S^3_{E_j}$ is in $\cH_1$ by Lemma \ref{021}, $S^3_E= \#_{j=1}^s S^3_{E_j}$ is in $\cH_1^\oplus$.

Then
\begin{align*} \hspace*{10pc}
w_I^{(n)}(D) & \ = \    I ([S^3,E]) \\
        &\ = \ I(S^3) - I(S^3_E) \\
        &\ = \ 0 \quad \text{because } S^3_E \in \cH_1^\oplus.
\end{align*}

It follows that $I$ is an invariant of degree $\le 2n-2$. Induction then shows that $I$ is an invariant of degree 0, or just a constant invariant. Hence $I(M)=I(S^3)=1$ for every integral homology
3-sphere $M$.

b) Suppose $I$ is a finite type invariant of degree $\le 2n$, and $D$ a Jacobi diagram of degree $n$. Let us restrict $I$ on the class $\cM(M)$. One has
\begin{align*} \hspace*{5pc}
w_I^{(n)}(D) & \ = \   I ([M\#S^3,E]) \\
&\ = \ I(M) - I(M\# S^3_E) \\
&\ = \ I(M) - I(M) \, I(S^3_E) \quad \text{because $I$ is multiplicative}\\
&\ =\ 0
\end{align*}
Hence again $I$ is an invariant of degree 0,
or $I$ is a constant invariant on every class $\cM(M)$.

Since $I(M)=1$ for every lens space of the form $L(p,1)$, it follows that if a rational homology sphere $M$ belongs to $\cM(N)$, where $N$ is the connected sum of a finite
number of lens spaces of the form $L(p,1)$, then $I(M)=1$.

Ohtsuki's lemma \cite{Ohtsuki_pi} says that for every rational homology sphere $M$, there are lens spaces $L(p_1,1), \dots, L(p_s,1)$ such that
the linking form of $N = M \# (  \#_{j=1}^s L(p_j,1))$ is the sum of the linking forms of a finite number lens spaces of the form $L(p,1)$. Since $I$ is multiplicative

$$I(N) = I(M) \prod_{j=1}^s I(L(p_j,1)).$$
With $I(N)=1= I(L(p_j,1))$, it follows that $I(M)=1$.
        \end{proof}

\subsubsection{Another proof of Theorem \ref{thm.rec} in the link case}

\begin{proof}[Proof of Theorem \ref{thm.rec} in the link case]
When $R$ is a field of characteristic 0,  the LMO invariant is universal for finite type. This fact can be reformulated as $W: \Gr_n \cA(\emptyset) \to
\cF_{2n}\cM(M)/\cF_{2n+1}\cM(M)$ is a bijection. This was proved for integral homology 3-spheres by Le \cite{Le2} and for general rational homology spheres by Habiro.
In particular, this result says that the part of degree $\le n$ of $\hLMO$ is a (universal) finite type invariant of degree $\le 2n$.

Note that $\hWg(\hLMO)$ and $\tau^g$ are multiplicative invariant with values in $\R[[\hbar]]$. By Proposition \ref{060} below, the part $\tau^\fg_{\le n}$ of degree less than or equal to $n$ of $\tau^\fg$ is a finite type invariant of degree $\le 2n$. Let $I= |H_1|^{\phi_+} \tau^g/\hWg(\hLMO)$. Then the part $I_{\le n}$ of degree less than or equal to $n$ is an invariant of degree less than or equal to $2n$. Clearly $I_{\le n}$ is multiplicative.
Moreover $I_{\le n}(M)=1$ if $M$ is obtained by surgery on knots by the knot case. Hence by Theorem \ref{022}, $I_{\le n}=1$. Since this holds true for every $n$ one has $I=1$, or
$\hWg(\hLMO) = \tau^\fg$.
\end{proof}

\section{Presentations of the perturbative invariants}
\label{sec.ppi}
In this section we discuss the perturbative invariants. In particular, we prove part (b) of Proposition \ref{n2003} and show that the degree $n$ part of
the perturbative is a finite type invariant of order $\le 2n$. We also give a informal way to explain how one can arrives at the formula of the perturbative invariant
given by Proposition \ref{n2003}.

\subsection{Perturbative expansion of a Gaussian integral}
\label{sec.pGint}

In this section, we explain how a Gaussian integral with a formal parameter in the exponent can be understood in perturbative expansions.
For the perturbative expansion of a Gaussian integral,
see also \cite[Appendix]{Aint1}.

\vskip 0.2pc
Suppose $V$ is a finite-dimensional Euclidean space, $f$ be a non-zero integer, $R\in S(V)= P(V^*)$.
The Gaussian integral

$$ I = \int_{V^*} e^{f \hbar |x|^2 /2}\, R(x)\, dx$$
does not make sense if $\hbar$ is a formal parameter. If $\hbar$ is a real number such that $f \hbar <0$, then the integral converges absolutely, and one can calculate the integral as follows.
A substitution $x = u /\sqrt{-2f\hbar}$ leads to
\begin{align*} \hspace*{5pc}
I & \ = \ \frac{1}{( -2 f \hbar )^{\dim V/2}} \,
 \int_{V^*} e^{-|u|^2 /4}\, R(\frac{u}{\sqrt{-2f\hbar}})\, du  \\
&\ = \ \left( \frac{2\pi}{-f\hbar}\right) ^{\dim V /2}\,  \eL^{(f)}_\fh(R) \quad \text{ by Proposition \ref{n02}}.
\end{align*}

If $\hbar$ is a formal parameter, the right-hand side still makes sense as an element in $\R[1/h]$. Thus we should declare

\begin{equation} \hspace*{5pc}
\int_{V^*} e^{f \hbar |x|^2 /2}\, R(x)\, dx
\ = \ \left( \frac{2\pi}{-f\hbar}\right) ^{\dim V /2}\,  \eL^{(f)}_\fh(R)
\label{eq.trans}
\end{equation}
for a formal parameter $\hbar$. Note that if $R \in S(V)[[\hbar]]$ then the right-hand side makes sense in $\R(\!(h)\!)$.

\subsection{Derivation of the perturbative invariants from the WRT invariant}
\label{sec.p_pi}

First we review the 3-manifold WRT invariant, for details see {\it e.g.} \cite{Le_isp}.
We again assume   $M$ is obtained by surgery on an algebraically split link $L$ with framing $\ff=(f_1,\dots,f_\ell)$. Let $L_0$ be the link $L$ with all framings 0, and $T$ is an algebraically string link (with 0 framing on each component) such that its closure is $L_0$.

For an $\ell$-tuple $(V_{\lambda_1-\rho}, \dots, V_{\lambda_\ell-\rho})$ of $\fg$-modules one can define the quantum %(Reshetikhin-Turaev)
link invariant
$Q^{\fg;V_{\lambda_1-\rho}, \dots, V_{\lambda_\ell-\rho}}(L_0)$ of the link $L_0$, (see \cite{RT}, we use here notations from the book \cite{Ohtsuki_book}). This invariant can be calculated though the Kontsevich invariant by results of \cite{Kassel,LM}:

\begin{equation} \hspace*{5pc}
Q^{\fg;V_{\lambda_1-\rho}, \dots, V_{\lambda_\ell-\rho}}(L_0)
\ = \ \big ( Z(T) \, \Delta^{(\ell)}(\nu)\big)(\lambda_1,\dots,\lambda_\ell).
\label{n3004}
\end{equation}
In particular, when $L_0=U$, the unknot with framing 0, $Q^{\fg;V_{\lambda_-\rho}}(U)$ is called the quantum dimension of $V_{\lambda_-\rho}$, denoted by $\qd(\lambda)$; its value is well-known:
\begin{equation}
\qd(\lambda)\ = \, \prod_{\al \in \Phi_+}
\frac{\big[ (\la,\al) \big]}{\big[ (\rho,\al) \big]},
\label{n6003}
\end{equation}
where $[n] := (q^{n/2}-q^{-n/2})/(q^{1/2}-q^{-1/2})$. Recall that we always have $q= e^\hbar$.

The quantum invariant of $L$ differs from that of $L_0$ by the framing factors, which will play the role of the exponential function in the Gaussian integral:
\begin{equation} \hspace*{4pc}
Q^{\fg;V_{\lambda_1-\rho}, \dots, V_{\lambda_\ell-\rho}}(L)
\ = \ \left(  \prod_{j=1}^\ell q^{f_j (|\lambda_j|^2 - |\rho|^2)/ 2} \right)\,
Q^{\fg;V_{\lambda_1-\rho}, \dots, V_{\lambda_\ell-\rho}}(L_0) .
\label{n3005}
\end{equation}
The normalization used in the definition of the WRT invariant is
$$
F_L(\lambda_1,\dots,\lambda_\ell) := \left( \prod_{j=1}^\ell \qd(\lambda_j) \right) Q^{\fg;V_{\lambda_1-\rho}, \dots, V_{\lambda_\ell-\rho}}(L).
$$
Using \eqref{n3004} and \eqref{n3005} one can show that
\begin{equation} \hspace*{4pc}
F_L(\lambda_1,\dots,\lambda_\ell)
\ = \ \Big( q^{-\sum_j f_j |\rho^2|/2}\Big)  \left( e^{\sum_j f_j \hbar |\lambda_j|^2 /2} \, R (\lambda_1,\dots,\lambda_\ell) \right),
\label{n3011}
\end{equation}
where $R= \od^{\otimes \ell}\, \Ups_\fg \left( \cQ^\fg(T)\right) = F_{L_0}$.

Suppose $q$ is a complex root of unity of order $r$.
Then it is known \cite{Le_isp} that
the function $F_L(\lambda_1,\dots,\lambda_\ell)$ is component-wise invariant under the
translation by $r\alpha$ for any $\alpha$ in the root lattice. Let $D_r \subset \fh^*$ be any fundamental domain of the translations by $r\alpha$ with $\alpha$ in the root lattice.
Then, with $q$ an $r$-th root of 1,
\begin{equation} \hspace*{10pc}
I(L) := \sum_{\lambda_j \in D_r}\, F_L(\lambda_1,\dots,\lambda_\ell)
\label{n3010}
\end{equation}
is invariant under the handle slide move. A standard normalization of $I(L)$
 gives us an invariant of 3-manifolds, which is the WRT invariant.

Because of the translational invariance of $F_L$, we could define the WRT invariant  if we replace $D_r$ by $N D_r$ in \eqref{n3010}, where $N$ is any positive integer. Let $N \to \infty$, we
should sum over all the weight lattice in \eqref{n3010} which does not converge. Instead, we use integral over $\fh^*$, {\it i.e.}, instead of $I(L)$ we consider the integral
 $$ \int_{(\fh^*)^\ell} F_L(\lambda_1,\dots,\lambda_\ell) d\lambda_1\dots d \lambda_\ell,$$
 which does not make sense. However, using $F_L(\lambda_1,\dots,\lambda_\ell)$ in \eqref{n3011}, the integral has the form of a Gaussian integral discussed in the previous section. According to
 \eqref{eq.trans}, the above integral should be a constant multiple of
following modification of $I(L)$:

$$ I_2(T,\ff) := \left( \prod_{j=1}^\ell  q^{-f_j |\rho|^2/2} \right) \, \eL^{(\ff)}_\fh \left(   \od^{\otimes \ell}\, \Upsilon_\fg ( \cQ^\fg(T)) \right),$$
which leads to the formula in Proposition \ref{n2003}.

\subsection{Proof of Proposition \ref{n2005}(b)}
\label{sec.p_Le}

First we review Le's formula of $\tau^\fg$,
for details, see \cite{Le_isp}.
As noted in the previous section,  as functions on $(\fh^*)$,
$$
F_{L_0} =  \od^{\otimes \ell} \,\Ups_\fg\big ( \cQ^\fg(T) \big).
$$
Let $\cO^{(f)}: P(\fh^*) = S(\fh) \to \R[1/\hbar]$ be the unique linear operator defined by
\begin{equation} \hspace*{3pc}
\cO^{(f)} (\beta^k)
\ = \ \begin{cases}
0 & \mbox{ if $k$ is odd, }
\\
\displaystyle{q^{- f |\rho|^2/2} (2d-1)!! \,
\big( - \frac{|\beta|^2}{f} \big)^d \hbar^{-d}}
& \mbox{ if $k = 2d$.}
\end{cases}
\label{n3012}
\end{equation}
for $\beta \in \fh$.
We also define its multi-linear extension
$$
\cO^{(\ff)}: S(\fh)^{\otimes \ell}[[\hbar]] \to \R(\!(h)\!), \quad \cO^{(\ff)} := \bigotimes_{j=1}^\ell \cO^{(f_j)}.
$$
Let
$$
I_2'(T,\ff) :=  \, \cO^{(\ff)} \left(   F_{L_0} \right)=  \, \cO^{(\ff)} \left(   \od^{\otimes \ell}\, \Upsilon_\fg ( \cQ^\fg(T)) \right).
$$
Then the $\tau^\fg(M)$ is given by \cite{Le_isp}
$$
\tau^\fg(M) =
\frac{I_2'(T,\ff)}{ \prod_{j=1}^\ell I_2'(\downarrow, \sn(f_j))} \, .
$$

To prove part (b) of Proposition \ref{n2003} one needs only to show that $I_2(T;\ff)= I_2'(T;\ff)$. It is enough to show that
\begin{equation} \hspace*{9pc}
\cO^{(f)}(g) \ = \ q^{- f |\rho|^2/2} \eL^{(f)}_\fg(g)
\label{n5002}
\end{equation}
 for every $g \in S(\fh)$. Since both operators $\cO^{(f)}$ and $\eL^{(f)}_\fg$ are linear and $G$-invariant, it is sufficient to consider the case when
$g= x_1^k$, where $x_1$ is the first vector of an orthonormal basis $x_1,\dots,x_n$ of $\fh$. In this case $\Delta_\fh= \sum \partial_{x_j}^2$, and one can easily calculate
$\eL^{(f)}_\fg(x_1^k) =
\exp( - \frac{\Delta_\fh}{-2f \hbar}) (x_1^k)|_{x_j=0}$,
\begin{align*} \hspace*{4pc}
\exp( \frac{\Delta_\fh}{-2f \hbar}) (x_1^k)|_{x_j=0}
&\ = \ \sum_d \frac{\Delta^d}{d! (-2f \hbar)^d} (x_1^k) \\
&\ = \ \begin{cases}
0 & \mbox{ if $k$ is odd, }
\\
\displaystyle{ (2d-1)!! \,
\big( - \frac{1}{f} \big)^d \hbar^{-d}}
& \mbox{ if $k = 2d$, }
\end{cases}
\end{align*}
which is precisely the right-hand side of \eqref{n3012} without the factor $q^{- f |\rho|^2/2}$ (with $\beta= x_1$). This proves \eqref{n5002}.
\qed

\subsection{The coefficients of $\tau^\fg$ are of finite type}

\begin{prop}\label{060}
The degree $n$ part of the perturbative invariant $\tau^\fg$
is a finite type invariant of degree $\le 2n$.
\end{prop}
\begin{rem} The proposition is a consequence of the main theorem. However, we used this proposition in the alternative proof of the main theorem in
Section \ref{sec.link_fti}. 
This is the reason why we give here a proof of the proposition independent of the main theorem.
\end{rem}
\begin{proof}
Let $M$ be a rational homology 3-sphere and
$E$ a collection of $2n+1$ disjoint $Y$-graphs in $M$.
We only need to prove that
$\tau^\fg([M,E])\in \hbar^{n+1} \BQ[[\hbar]]$.

By taking connecting sum with lens spaces, we assume that the pair $(M,E)$ can be obtained from $(S^3,E)$ by surgery along an algebraically split link $L\subset S^3$.
By adding trivial knots with framing $\pm1$ (which are unlinked with $L$) to $L$ if needed, we can assume that the leaves of $E\in S^3$ form a zero-framing trivial link. Let $L_0$ be the link $L$ with 0 framing, and choose a string link $T$ in  a cube such that $L_0$ is the closure of $T$. We can assume that $E$ is also in the cube.

For a sub-collection  $E'\subset E$ let $L_{E'}$ be the link obtained by surgery of $S^3$ along $E'$ (see \cite{Ha2, GGP}). We define similarly $(L_0)_{E'}$ and $T_{E'}$. Clearly $(L_0)_{E'}$ is the closure of $T_{E'}$.

For every link $L$ and every $Y$-graph $C$  whose leaves are a zero-framing trivial link,
the moves from $L$ to $L_C$ is a repetition of the Borromeo move
(see \cite{Ha2,GGP}):

\begin{figure}[htb]
$$ \psdraw{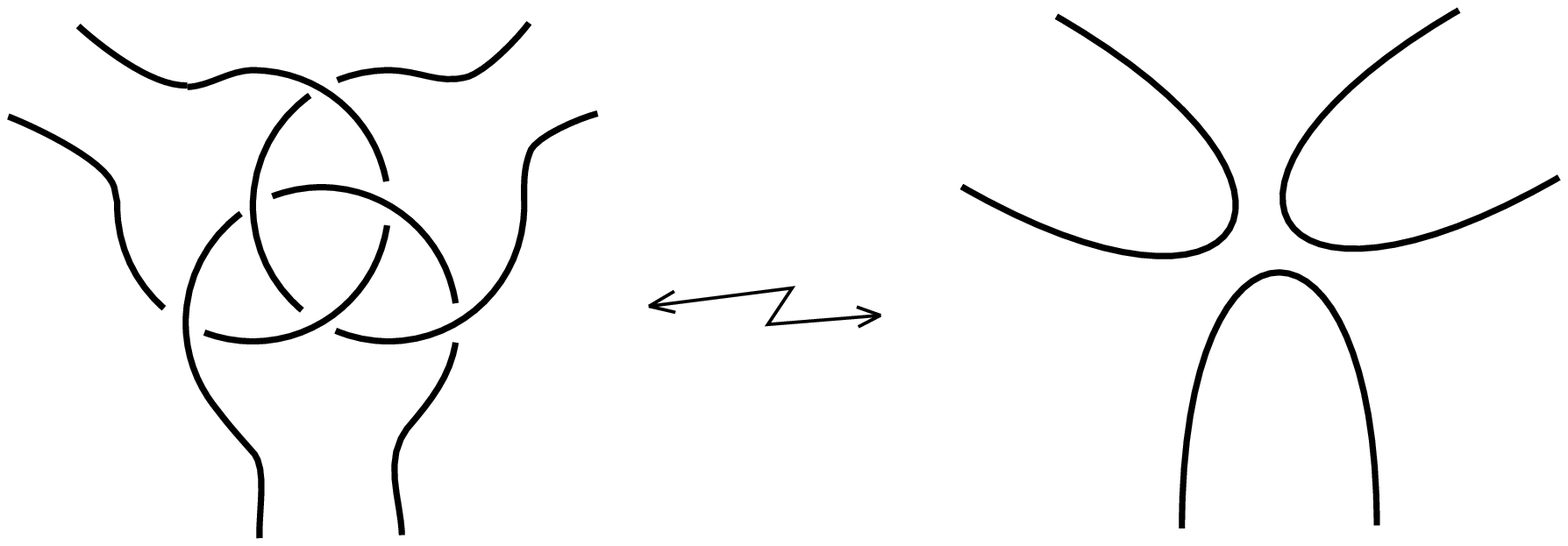}{2.6in} $$
\caption{\label{move23}}
\end{figure}

Hence, by \cite[Lemma 5.3]{Le2}, $Z(T-T_C)$ has i-degree $ \ge  1$.
Here $x \in \cA(\sqcup  ^\ell \downarrow)$ has i-degree $\ge k$ if it is a linear combinations of Jacobi diagrams with at least $k$ trivalent vertices.
 It follows that $\cZ([T,E])$ has i-degree $\ge 2n+1$, where $[T,E] := \sum_{E'\subset E} (-1)^{|E'|} T_{E'}$.

 Note that all the links $L_{E'}, E'\subset E$ are algebraically split, having the same number of components, and having same the framings $\ff=(f_1,\dots,f_\ell)$.
 By definition, one has
\begin{align*}  \hspace*{10pc}
[M,E]
& \ = \  \sum_{E'\subset E} (-1)^{|E'|} (S^3)_{L_{E'}}
\end{align*}

Hence
\begin{align}  \hspace*{5pc}
\tau^\fg([M,E]) \ & = \
\sum_{E'\subset E} (-1)^{|E'|} \tau^\fg \left ((S^3)_{L_{E'}} \right) \notag\\
&= \
\sum_{E'\subset E} (-1)^{|E'|}  \frac{I_2(T_{E'},\ff)}{ \prod_{j=1}^\ell I_2\big(\downarrow,\sn(f_j)\big)}  \notag\\
& = \, \frac{I_2([T,E],\ff)}{ \prod_{j=1}^\ell I_2\big(\downarrow,\sn(f_j)\big)}  \notag\\
&=
\frac{ \left( \prod_{j=1}^\ell  q^{-f_j |\rho|^2/2} \right) \, \eL^{(\ff)}_\fh \left(   \od^{\otimes \ell}\, \cZ([T,E]) \right) }
{ \prod_{j=1}^\ell I_2\big(\downarrow,\sn(f_j)\big)} \, .
\label{n6000}
\end{align}

 By Lemma \ref{n6001} below, since $\cZ([T,E])$ has i-degree $\ge 2n+1$, the numerator of \eqref{n6000} belongs to $\hbar^{-\ell \phi_+ + n+1} \R[[\hbar]]$, while the denominator has the form
 $\hbar^{-\ell \phi_+} u$, where $u$ is a unit in $\R[[h]]$. It follows that the right hand side of \eqref{n6000} belongs to $\hbar^{n+1} \R[[\hbar]]$. \end{proof}

 \begin{lem} a) Suppose $D\in \cA(\sqcup^\ell \downarrow)$ is a Jacobi diagram having $\ge 2n+1$ trivalent vertices, then $\eL^{(\ff)}_\fh \left(   \od^{\otimes \ell}\, \hWg (D)) \right)\in \hbar^{-\ell \phi_+ +n+1} \R[[\hbar]]$.

 b) The lowest degree of $\hbar$ in $ I_2\big(\downarrow,\pm 1\big)\in \R[[h]]$ is  $\hbar^{-\phi_+} $, i.e. $ \hbar^{\phi_+} I_2\big(\downarrow,\pm 1\big)$ is invertible in $\R[[\hbar]]$.
 \label{n6001}
 \end{lem}
 \begin{proof} a) Suppose $D$ has degree $d$. Then $D$ has $2d$ vertices, among which  $2d-2n-1$ are univalent. It follows that $\Wg(D)$, as element of $U(\fg)^{\otimes \ell}$, has degree  $\le (2d-2n-1)$, and, as a function on $(\fh^*)^\ell$, is a polynomial
 of degree $\le \ell \phi_+ + (2d-2n-1)$, see \cite{Le_isp}. Hence the degree of $\od^{\otimes \ell } \, \Wg(D)$ is $ \le 2\ell \phi_+ + 2d-2n-1$.  Recall that $\eL^{(\ff)}_\fh(g)$ lower the degree of $\hbar$ by at most half the degree of $g$. 
The degree of $\hbar$ in 
$\eL^{(\ff)}_\fh \big( \od^{\otimes \ell}\, \hWg (D) \big)
= \hbar^d \, \eL^{(\ff)}_\fh \big( \od^{\otimes \ell}\, \Wg (D) \big)$ 
is at least $d - \frac{1}{2}(2\ell \phi_+ + 2d-2n-1) = 1/2 + n - \ell \phi_+$.
Hence $\eL^{(\ff)}_\fh \big( \od^{\otimes \ell}\, \hWg (D) \big) 
\in \hbar^{-\ell \phi_+ +n+1} \R[[\hbar]]$.

 b) By definition

 $$ I_2\big(\downarrow,\pm 1\big) = q^{\mp |\rho|^2/2}\, \eL^{(\pm 1})_\fh (\od \, \cZ(\downarrow)).$$
 For the trivial knot everything can be calculated explicitly. One has $\od \, \cZ(\downarrow)= (\qd)^2$, and using \eqref{n6003} one can easily show that
 $$ \od \, \cZ(\downarrow) = \od^2 \left( 1 + \sum_{k=1}^\infty g_k \hbar^{2k}\right),$$
 where $g_k$ has degree exactly $2k$. Thus
$$
q^{\pm  |\rho|^2/2}\, I_2\big(\downarrow,\pm 1\big)
= \eL^{(\pm 1)}_\fh  (\od)^2 +  \eL^{(\pm 1)}_\fh
\left ( \sum_{k=1}^\infty g_k \od^2 \hbar^{2k} \right).
$$
 since $\deg(g_k)=2k$, and $\deg(\od^2)= 2\phi_+$, the second term belongs to $ \hbar^{1-\phi_+} \R[[\hbar]]$, while the first term is
$$ 
\eL^{(\pm 1)}_\fh  (\od)^2 
\, = \, \hbar^{-\phi_+}
\frac{\Delta_\fh^{\phi_+} (\od^2)}{(\phi_+)! (\mp 2)^{\phi_+}}
\, = \, \hbar^{-\phi_+} \frac{ c_\fg}{(\mp 2)^{\phi_+}} \, .
$$
 Since $c_\fg\neq 0$ and $q^{\pm  |\rho|^2/2}$ is invertible, we conclude that $ h^{\phi_+} I_2\big(\downarrow,\pm 1\big)$ is invertible in $\R[[\hbar]]$.
 \end{proof}

\normalsize

Research Institute for Mathematical Sciences,
Kyoto University, Sakyo-ku, Kyoto, 606-8502, Japan

E-mail address: {\tt marron@kurims.kyoto-u.ac.jp}

\vskip 1pc

School of Mathematics, 686 Cherry Street, Georgia Tech, Atlanta,
GA 30332, USA

E-mail address: {\tt letu@math.gatech.edu}

\vskip 1pc

Research Institute for Mathematical Sciences,
Kyoto University, Sakyo-ku, Kyoto, 606-8502, Japan

E-mail address: {\tt tomotada@kurims.kyoto-u.ac.jp}


\begin{thebibliography}{99}

\normalsize

\ifx\undefined\bysame
\newcommand{\bysame}{\leavevmode\hbox to3em{\hrulefill}\,}
\fi

\bibitem{BarNatan}
Bar-Natan, D.,
{\it On the Vassiliev knot invariants},
Topology {\bf 34} (1995) 423--472.

\bibitem{BGRT_w}
Bar-Natan, D., Garoufalidis, S., Rozansky, L., Thurston, D.P.,
{\it Wheels, wheeling, and the Kontsevich integral of the unknot},
Israel J. Math. {\bf 119} (2000) 217--237.

\bibitem{Aint1}
\bysame,
{\it The Aarhus integral of rational homology 3-spheres I:
A highly non trivial flat connection on $S^3$},
Selecta Math. (N.S.) {\bf 8} (2002) 315--339.

\bibitem{Aint2}
\bysame,
{\it The Aarhus integral of rational homology 3-spheres II:
Invariance and universality},
Selecta Math. (N.S.) {\bf 8} (2002) 341--371.

\bibitem{Aint3}
\bysame,
{\it The Aarhus integral of rational homology 3-spheres III:
The relation with the Le-Murakami-Ohtsuki invariant},
Selecta Math. (N.S.) {\bf 10} (2004) 305--324.

\bibitem{BaLa}
Bar-Natan, D., Lawrence, R.,
{\it A rational surgery formula for the LMO Invariant},
Israel J. Math. {\bf 140} (2004) 29--60.

\bibitem{BGV}
Berline, N.,  Getzler, E.,  Vergne, M., {\it  Heat kernels and Dirac operators}.   Springer-Verlag, Berlin, 2004.
\bibitem{BLT}
Bar-Natan, D., Le, T.T.Q., Thurston, D.P.,
{\it Two applications of elementary knot theory to Lie algebras and
Vassiliev invariants},
Geometry and Topology {\bf 7} (2003) 1--31.

\bibitem{BBL}
Beliakova, A., Buehler, I., Le, T.,
{\it A unified quantum $SO(3)$ invariant for rational homology 3-spheres},
arXiv:0801.3893.

% \bibitem{CP} Chari, V., Pressley, A., {\it A guide to quantum groups}, Cambridge University Press, Cambridge, 1994.

\bibitem{FdV} Freudenthal, H., de Vries, H., {\it Linear Lie groups}, Pure and Applied Mathematics {\bf 35}, Academic Press, New York-London 1969.

\bibitem{GGP} Garoufalidis, S., Goussarov, M.,  Polyak., {\it  Calculus of clovers and finite type invariants of 3-manifolds},   Geom. Topol.  {\bf 5}  (2001), 75--108 (electronic).
\bibitem{Ha2} Habiro K., {\it Claspers and finite type invariants of links}, Geom. Topol. {\bf 4} (2000) 1–83.
\bibitem{Habiro_q}
Habiro, K.,
{\it On the quantum $sl_2$ invariants of knots and integral homology spheres},
Invariants of knots and 3-manifolds (Kyoto 2001), 161--181,
Geom. Topol. Monogr. {\bf 4}, Geom. Topol. Publ., Coventry, 2002.

\bibitem{Habiro_WRT}
\bysame,
{\it A unified Witten-Reshetikhin-Turaev invariant for integral homology
spheres}, Invent. Math. {\bf 171} (2008) 1--81.

\bibitem{HabiroLe}
Habiro, K., Le, T.T.Q.,
in preparation.

\bibitem{Hall} Hall, B.,  {\it Lie groups, Lie algebras, and representations: an elementary introduction}, Graduate Text in Mathematics {\bf 222}. Springer--Verlag, New York, 2003

%\bibitem{HC}
%Harish-Chandra,
%{\it Fourier transforms on a semisimple Lie algebra. I},
%Amer. J. Math. {\bf 79} (1957) 193--257.

\bibitem{Helgason}
Helgason, S.,
{\it Groups and geometric analysis.
Integral geometry, invariant differential operators, and spherical functions},
Mathematical Surveys and Monographs {\bf 83}.
American Mathematical Society, Providence, RI, 2000.

\bibitem{HT}
Howe, R., Tan, E.-C.,
{\it Non-abelian harmonic analysis, Applications of ${\rm SL}(2,R)$},
Universitext. Springer-Verlag, New York, 1992.

\bibitem{Kassel}
Kassel, C.,
{\it Quantum groups},
Graduate Texts in Mathematics {\bf 155}. Springer-Verlag, New York, 1995.

\bibitem{Kontsevich}
Kontsevich, M.,
{\it Vassiliev's knot invariants},
Adv. in Sov. Math {\bf 16(2)} (1993) 137--150.

\bibitem{Kuriya}
Kuriya, T.,
{\it On the LMO conjecture},
preprint, arXiv:0803.1732.

\bibitem{Le2}
Le, T.T.Q.,
{\it An invariant of integral homology 3-spheres which is universal
for all finite type invariants},
AMS translation series 2, Eds. V. Buchtaber and S. Novikov {\bf 179}
(1997) 75--100.

\bibitem{Le_isp}
\bysame,
{\it Quantum invariants of 3-manifolds: integrality, splitting, and
perturbative expansion},
Proceedings of the Pacific Institute for the Mathematical Sciences Workshop
``Invariants of Three-Manifolds'' (Calgary, AB, 1999).
Topology Appl. {\bf 127} (2003) 125--152.

\bibitem{LM}
Le, T.T.Q., Murakami, J.,
{\it The universal Vassiliev-Kontsevich invariant for framed oriented links},
Compositio Math. {\bf 102} (1996) 41--64.

\bibitem{LMO}
Le, T.T.Q., Murakami, J., Ohtsuki, T.,
{\it On a universal perturbative invariant of 3-manifolds},
Topology {\bf 37} (1998) 539--574.

\bibitem{Matveev}
Matveev, S.,
{\em Generalized surgery of 3-dimensional manifolds and
representations of homology 3-spheres} (in Russian), Mat. Zametki,
{\bf 42} (1987), 268--275.

\bibitem{Ohtsuki_pi}
Ohtsuki, T.,
{\it A polynomial invariant of rational homology 3-spheres},
Invent. Math. {\bf 123} (1996) 241--257.

\bibitem{Ohtsuki_rec}
\bysame,
{\it The perturbative $SO(3)$ invariant of rational homology 3-spheres
recovers from the universal perturbative invariant},
Topology {\bf 39} (2000) 1103--1135.

\bibitem{Ohtsuki_book}
\bysame,
{\it Quantum invariants, --- A study of knots, 3-manifolds, and their sets},
Series on Knots and Everything {\bf 29}.
World Scientific Publishing Co., Inc., 2002.

\bibitem{RT1} N. Yu. Reshetikhin, V. Turaev,
{\it Ribbon graphs and their invariants derived from quantum
groups}, Commun. Math. Phys., {\bf 127} (1990), 1--26.
\bibitem{RT}
Reshetikhin, N., Turaev, V.G.,
{\it Invariants of $3$-manifolds via link polynomials and quantum groups},
Invent. Math. {\bf 103} (1991) 547--597.

\bibitem{Rozansky}
Rozansky, L.
{\it Witten's invariants of rational homology spheres at prime values of $K$
and trivial connection contribution},
Comm. Math. Phys. {\bf 180} (1996) 297--324.

% \bibitem{Wakimoto} Wakimoto, M., {\it Infinite-dimensional Lie algebras}, Translations of Mathematical Monographs {\bf 195}. Iwanami Series in Modern Mathematics. American Mathematical Society, Providence, RI, 2001.

\bibitem{Witten}
Witten, E.,
{\it Quantum field theory and the Jones polynomial},
Comm. Math. Phys. {\bf 121} (1989) 351--399.

\end{thebibliography}
\end{document}